\documentclass[11pt,a4paper]{amsart}
\usepackage{amsmath,amsfonts,amssymb,amscd,amsthm}
\usepackage{graphicx}
\usepackage{color}
\usepackage[francais]{babel}
\newtheorem{theoreme}[equation]{Th\'eor\`eme}
\newtheorem{lemme}[equation]{Lemme}
\newtheorem{proposition}[equation]{Proposition}
\newtheorem{definition}[equation]{D\'efinition}
\newtheorem{corollaire}[equation]{Corollaire}
\newtheorem{conjecture}[equation]{Conjecture}
\newtheorem{notation}[equation]{Notation}
\theoremstyle{remark}
\newtheorem{remarque}[equation]{Remarque}
\numberwithin{equation}{section}

\DeclareMathOperator{\ppcm}{ppcm}
\DeclareMathOperator{\pgcd}{pgcd}
\newcommand{\bbC}{{\mathbb C}}
\newcommand{\bbZ}{{\mathbb Z}}
\newcommand{\Sgot}{{\mathfrak S}}
\newcommand{\bs}{{\bf s}}
\newcommand{\CO}{{\mathcal O}}
\newcommand{\uc}{{\underline c}}
\newcommand{\ur}{{\underline r}}
\newcommand{\us}{{\underline s}}
\newcommand{\ut}{{\underline t}}
\newcommand{\uw}{{\underline w}}
\newcommand{\ux}{{\underline x}}
\newcommand{\uP}{{\underline P}}
\newcommand{\inv}{^{-1}}
\newcommand{\cf}{{\it cf.}}
\newcommand{\ie}{{\it i.e.}}
\def\etc{{\it etc\dots}}
\newcommand{\Bn}{B(\tilde A_{n-1})}
\newcommand{\quasiparab}{quasi-parabolique}
\begin{document}
\author{F.~Digne}
\title[Pr\'esentations duales des groupes de tresses de type affine]
{Pr\'esentations duales des groupes de tresses de type affine $\tilde A$}
\address{LAMFA, Universit\'e de Picardie-Jules Verne
33, Rue Saint-Leu 80039 Amiens Cedex France}
\email{digne@u-picardie.fr}
\urladdr{http://www.lamfa.u-picardie.fr/digne}
\subjclass[2000]{Primary 20F36; Secondary 20F05}
\keywords{Garside monoid, dual monoid, non-crossing partitions, affine
braids}

\maketitle
\section*{Abstract}
Artin-Tits groups of spherical type have two well-known Garside structures,
coming respectively from the divisibility properties of the classical
Artin monoid and of the dual monoid. For general Artin-Tits groups, the
classical monoids have no such Garside property. In the present paper we
define dual monoids for all Artin-Tits groups and we prove that for the type
$\tilde A_n$ we get a (quasi)-Garside structure.  Such a structure provides
normal forms for the Artin-Tits group  elements and allows to solve some
questions such as to determine the centralizer of a power of the Coxeter
element in the Artin-Tits group.

More precisely, if $W$ is a Coxeter group, one can consider the length $l_R$
on $W$ with respect to the generating set $R$ consisting of all reflections.
Let $c$ be a Coxeter element in $W$ and let $P_c$ be the set of elements $p\in
W$ such that $c$ can be written $c=pp'$ with $l_R(c)=l_R(p)+l_R(p')$.  We
define the monoid $M(P_c)$ to be the monoid generated by a set $\underline
P_c$ in one-to-one correspondence, $p\mapsto \underline p$, with $P_c$ with
only relations $\underline{pp'}=\underline p.\underline p'$ whenever $p$, $p'$
and $pp'$ are in $P_c$ and $l_R(pp')=l_R(p)+l_R(p')$.  We conjecture that the
group of quotients of $M(P_c)$ is the Artin-Tits group associated to $W$ and
that it has a simple presentation (see \ref{conjecture} (ii)). These
conjectures are known to be true for spherical type Artin-Tits groups. Here we
prove them for Artin-Tits groups of type $\tilde A$. Moreover, we show that
for exactly one choice
of the Coxeter element (up to diagram automorphism)
we obtain a (quasi-) Garside monoid.
The proof makes use of non-crossing paths in an annulus which are the
counterpart in this context of the non-crossing partitions used for type $A$.

\section{Introduction}
Soit $(W,S)$ un syst\`eme de Coxeter quelconque (avec $S$ fini). 
Soit $R$ l'ensemble des r\'eflexions de $W$ (c'est-\`a-dire des conjugu\'es des
\'el\'ements de $S$).
Appelons ``longueur de r\'eflexion'' $l_R(w)$ d'un \'el\'ement $w\in W$
le nombre minimum de termes dans une d\'ecomposition de $w$ en
produits de r\'eflexions. Nous dirons que $v\in W$ divise $w\in W$,
not\'e $v\preccurlyeq w$,  si $w=vv'$ avec $l_R(w)=l_R(v)+l_R(v')$.
Comme $l_R$ est invariant par conjugaison ceci est \'equivalent \`a $w=v''v$ avec
$l_R(w)=l_R(v'')+l_R(v)$, autrement dit il n'y a pas lieu de distinguer entre
diviseurs \`a gauche et \`a droite. La relation $\preccurlyeq$ est clairement une
relation d'ordre.
Fixons un \'el\'ement $c$ de Coxeter, c'est-\`a-dire un
produit de tous les \'el\'ements de $S$ dans
un certain ordre; consid\'erons l'ensemble $P_c$ des diviseurs de $c$ et
$M(P_c)$  le mono{\"\i}de engendr\'e par un ensemble $\underline P_c=\{\underline
w\,\mid\,w\in P_c\}$ en bijection avec
$P_c$ avec comme relations $\underline w.\underline w'=\underline{ww'}$ si $w$ et
$w'$ sont des \'el\'ements de $P_c$ tels que $ww'\in P_c$ et $l_R(ww')=l_R(w)+l_R(w')$ (\cf\ 
\cite[section 2]{BDM} et \cite[0.2 et 0.4]{B}).
La longueur $l_R$ s'\'etend \`a $M(P_c)$ en une longueur additive. Les notions de
divisibilit\'e \`a gauche  et \`a droite
dans le mono{\"\i}de \'etendent la relation de divisibilit\'e de $P_c$.
Nous conjecturons
\begin{conjecture}{\label{conjecture}}
\begin{enumerate}
\item Le groupe des fractions de $M(P_c)$ est isomorphe au groupe des
tresses d'Artin-Tits associ\'e \`a $W$.
\item Le mono{\"\i}de $M(P_c)$ (resp. son groupe de fractions)
a la pr\'esentation suivante comme mono{\"\i}de (resp. comme groupe):
l'ensemble des g\'en\'e\-rateurs est
$\{\underline r\mid r\in R\cap P_c\}$ et pour chaque
couple $(r,t)\in\ R^2$ tel que $rt$ divise $c$ on a la relation
$\ur.\ut=\underline{rtr}.\ur$.

\end{enumerate}
\end{conjecture}
Cette conjecture a \'et\'e prouv\'ee pour les groupes de Coxeter de type $A_n$ dans
\cite{BiKoLee} et \cite{BDM} et pour les autres types sph\'eriques dans \cite{B}.
Dans cet article nous allons prouver cette conjecture dans le cas o\`u $W$ est
de type $\tilde A_n$. Nous montrons de plus que
pour exactement un choix $c_0$ de  l'\'el\'ement de Coxeter (\`a
automorphisme du diagramme pr\`es)
le mono{\"\i}de obtenu est un treillis pour la divisibilit\'e (structure
quasi-Garside), ce qui implique l'existence de formes normales dans son groupe de fractions
et donne ainsi entre autres cons\'equences
une nouvelle solution au probl\`eme des mots dans les groupes
d'Artin-Tits de type $\tilde A_n$. Les preuves utilisent en particulier
l'introduction d'objets ``sans croisements'' dans une couronne,
et leur interpr\'etation comme tresses dans un cylindre.

Les sym\'etries des pr\'esentations obtenues d\'efinissent des automorphismes du groupe d'Artin-Tits dont on peut calculer les points fixes gr\^ace \`a la
th\'eorie des structures de Garside. Cela permet par exemple d'obtenir le centralisateur d'une puissance de l'\'el\'ement de Coxeter $c_0$,
qui s'av\`ere \^etre un groupe d'Artin-Tits de type $B$.

Dans la section 2 nous \'etudions la longueur $l_R$ et d\'eterminons
$P_c$ pour un
groupe de Coxeter de type $\tilde A$. Nous en donnons une interpr\'etation
topologique par des chemins sans croisement dans une couronne.
Dans la section 3 nous donnons une
pr\'esentation du mono{\"\i}de dual $M(P_c)$. Dans la section 4 nous prouvons la
conjecture \ref{conjecture} (i) dans le cas $\tilde A$. Dans la section 5 nous
prouvons qu'on a une structure de Garside et en donnons un certain nombre de
cons\'equences.

\section{Longueur de r\'eflexion; diviseurs d'un \'el\'ement de Coxeter} 
On sait par Dyer \cite{Dyer} que la longueur de r\'eflexion d'un \'el\'ement $w$
d'un groupe de Coxeter quelconque est \'egale au nombre minimum
de termes qu'il faut effacer dans une suite minimale d'\'el\'ements de $S$ de
produit \'egal \`a $w$ pour que le produit des termes restants soit \'egal \`a 1.
C'est aussi la longueur d'un plus court chemin (croissant)
de 1 \`a $w$ dans le graphe de Bruhat du groupe de Coxeter.
Ceci montre imm\'ediatement
par exemple que la longueur de r\'eflexion d'un \'el\'ement de Coxeter est \'egal \`a sa
longueur de Coxeter. Cela montre aussi que la restriction \`a un sous-groupe parabolique
de $l_R$ est \'egale \`a la longueur de r\'eflexion dans ce sous-groupe.
Soit $W$ un groupe de Coxeter de type $\tilde A_{n-1}$.
Pour \'etudier $M(P_c)$ nous avons besoin d'une formule explicite pour $l_R$ qui
nous permettra en particulier de 
d\'eterminer quels sont les diviseurs d'un \'el\'ement de Coxeter. Pour
cela nous utilisons une repr\'esentation de $W$ comme sous-groupe du groupe des permutations
de $\bbZ$. Pour expliquer cette interpr\'etation nous introduisons les d\'efinitions
suivantes:
\begin{definition}
\begin{itemize}
\item[$\bullet$] Nous dirons
qu'une permutation de $\bbZ$ est $n$-p\'eriodique si elle v\'erifie
$w(x+n)=w(x)+n$ pour tout $x\in \bbZ$.
\item[$\bullet$] Si $w$ est une permutation $n$-p\'eriodique
de $\bbZ$, nous appelons
d\'ecalage de $w$ l'entier $\dfrac1n\sum_{x=1}^{x=n}(w(x)-x)$.
\end{itemize}
\end{definition}
On sait (\cf\  par exemple \cite{shi}) que $W$ s'identifie aux permutations
$n$-p\'eriodiques de $\bbZ$ de d\'ecalage nul.
On peut le voir de la fa\c con suivante: le groupe de Coxeter
affine de type $\tilde A_{n-1}$ est le produit
semi-direct du r\'eseau des racines de type $A_{n-1}$
par le groupe de Coxeter de ce type.
Autrement dit $W=\Sgot_n\ltimes\{(a_i)\in\bbZ^n\mid\sum a_i=0\}$.
Notons $\overline w$ l'image dans $\Sgot_n$ de $w\in \Sgot_n\ltimes \bbZ^n$.
On peut identifier le groupe
$\Sgot_n\ltimes\bbZ^n$ au groupe des permutations $n$-p\'eriodiques de $\bbZ$,
l'image de $k\in\{1,\ldots n\}$ par $w=\overline w.(a_1,\ldots,a_n)$ \'etant
$\overline w(k)+na_k$. 
Alors le d\'ecalage d'une permutation $n$-p\'eriodique est
donn\'e par le morphisme
$\overline w.(a_1,\ldots,a_n)\mapsto a_1+\cdots+a_n$ et
$W$ s'identifie au noyau de ce morphisme.

Si $w$ est une
permutation $n$-p\'eriodique de $\bbZ$, 
l'ensemble des orbites de $w$ est invariant par translation
de $n$. Rappelons qu'on appelle support d'une
permutation $w$ le compl\'ementaire dans
$\bbZ$ de l'ensemble des points fixes de
$w$.
\begin{definition}\label{decalage}
Soit $w$ une permutation $n$-p\'eriodique de $\bbZ$;
nous appelons d\'ecalage d'une orbite $\CO$ de $w$ l'entier
$\dfrac1n\sum_{x\in\CO\cap\{1,\ldots,n\}}(w(x)-x)$.
\end{definition}
Le d\'ecalage de $\CO$ est donc l'entier $h$ tel que pour tout $a\in\CO$
on ait $w^k(a)=a+hn$ o\`u $k$ est le cardinal de l'image de
l'orbite modulo $n$.
Une orbite est finie si et seulement si elle est de d\'ecalage nul. 
Le d\'ecalage d'une permutation
$n$-p\'eriodique est la somme des d\'ecalages de ses orbites, donc
une permutation $n$-p\'eriodique de $\bbZ$
est dans $W$ si et seulement si la somme des
d\'ecalages des orbites infinies est nulle.

On peut d\'ecomposer un \'el\'ement de $W$ en produit d'\'el\'ements de
$W$ ayant deux \`a deux des supports disjoints.
On s'int\'eresse aux
d\'ecompositions maximales de cette forme.
Ceci nous am\`ene \`a poser la d\'efinition suivante:
\begin{definition}\label{cycle}
\begin{itemize}
\item[(i)] Nous dirons qu'une permutation $n$-p\'eriodique est un cycle si elle
n'a qu'une orbite non triviale \`a translation de $n$ pr\`es.
\item[(ii)] Nous dirons qu'un \'el\'ement de $W$ diff\'erent de l'identit\'e
est un pseudo-cycle si la restriction de $w$ \`a
toute partie stricte stable de son support n'est pas dans $W$.
\end{itemize}
\end{definition}
Toute permutation p\'eriodique est produit de cycles disjoints et cette
d\'ecomposition est unique \`a l'ordre pr\`es. Mais un cycle n'est dans $W$ que si
son d\'ecalage est nul, c'est-\`a-dire si les orbites de ce cycle sont finies. 
Tout \'el\'ement de $W$ est
produit de pseudo-cycles de supports disjoints, mais il n'y a pas unicit\'e des pseudo-cycles
ayant des orbites infinies (voir par exemple \ref{non-additivite}
ci-dessous). Un pseudo-cycle ou bien est un cycle,
ou bien a toutes
ses orbites non triviales infinies, de d\'ecalage total nul,
et dans ce cas toute sous-famille de ses orbites a un d\'ecalage total non nul.
\begin{notation}\label{(a,b,c)}

\begin{itemize}
\item[$\bullet$]
Un cycle sera repr\'esent\'e
sous la forme
$(a,b,c,\ldots,l)_{[h]}$ o\`u $h$ est le d\'ecalage du cycle. Cette notation
signifiant que
$a$, $b$, $c$,
\dots, $l$ sont tous distincts modulo $n$ et que
$w$ envoie $a$ sur $b$, $b$ sur $c$, \dots\ et
$l$ sur $a+hn$.
Pour simplifier les notations, si le d\'ecalage est nul, nous
omettrons l'indice $[0]$.
\item[$\bullet$]
Nous repr\'esenterons  toute permutation $n$-p\'eriodique de $\bbZ$
comme un produit de cycles de supports disjoints non vides.
\item[$\bullet$]
On pose $s_i=(i,i+1)$ pour $i=1,\ldots, n$. Les $s_i$ sont les g\'en\'erateurs de
Coxeter de $W$.
\end{itemize}
\end{notation}
Avec ces conventions
les r\'eflexions sont les \'el\'ements $(a,b)$ avec $a$ et $b$ quelconques distincts modulo $n$.
On obtient exactement une fois chaque r\'eflexion si on impose de plus $a<b$ et $a\in\{1,\ldots,n\}$. 

Nous utiliserons \`a plusieurs reprises
la formule suivante qui r\'esulte d'un calcul imm\'ediat.
\begin{lemme}\label{calcul}
On d\'esigne par $a_1,\ldots,a_l$ des entiers distincts
modulo $n$ et par $h$ et $k$ des entiers quelconques.
On a
\begin{multline*}
(a_1,\ldots,a_i,\ldots,a_l)_{[h]}(a_1,a_i+kn)=\\
(a_{i+1},\ldots,a_l,a_1+hn)_{[h+k]}(a_2,a_3,\ldots,a_i)_{[-k]}
\end{multline*}
\end{lemme}
\begin{definition}
Pour toute partie $T\subset\bbZ$ invariante par translation de $n$, notons
$\tilde W_T$ le fixateur dans $W$ du compl\'ementaire de $T$ et $W_T$ le
sous-groupe des \'el\'ements de $\tilde W_T$ n'ayant pas d'orbite infinie.
Un sous-groupe du type $W_T$ ou $\tilde W_T$ sera dit ``\quasiparab\ ''.
\end{definition}
Remarquons que $\tilde W_T$ et $W_T$ sont des groupes de Coxeter
de types respectifs
$\tilde A_{|\overline T|-1}$ et $A_{|\overline T|-1}$ si $\overline T$ est
l'image de $T$ modulo $n$.
\begin{notation}\label{nu et kappa}
Pour $w\in W$ notons
$\nu(w)$ le nombre de classes modulo $n$ d'orbites de $w$
et notons $\kappa(w)$ le plus grand entier $k$ tel
que $w=w_1w_2\ldots w_k$ o\`u
les $w_i\in W$ sont (des pseudo-cycles) de supports disjoints.
\end{notation}
Remarquons que $\kappa(w)$ est le plus grand entier $k$ tel que
$w$ soit dans le produit direct
de $k$ sous-groupes \quasiparab  s. On a $\kappa(w)\leq\nu(w)$.

La proposition suivante donne une formule explicite pour la longueur de
r\'eflexion dans $W$.
\begin{proposition}\label{l_R} Pour $w\in W$ on a
 $l_R(w)=n+\nu(w)-2\kappa(w)$.
\end{proposition}
\begin{proof}
On pose $f(w)=n+\nu(w)-2\kappa(w)$. On montre les trois propri\'et\'es
suivantes:
\begin{itemize}
\item[(1)] Si $f(w)\geq 1$, il existe $r\in R$ tel que $f(wr)=f(w)-1$.
\item[(2)] Si $f(w)=0$ alors $w=1$.
\item[(3)] $f(w)\leq l_R(w)$.
\end{itemize}
Les propri\'et\'es (1) et (2) donnent l'in\'egalit\'e $l_R(w)\leq f(w)$, d'o\`u la
proposition.

Prouvons (1).
Si $w$ a un cycle de la forme $(a_1,a_2,\ldots,a_l)_{[h]}$ avec $l\neq 1$,
la multiplication par $r=(a_1,a_2)$, augmente $\nu(w)$ et $\kappa(w)$ de 1
(\cf\  \ref{calcul}).
Si tous les cycles de $w$ sont de la forme
$(a_1)_{[h]}$, soit $w_i$  un des facteurs diff\'erents de l'identit\'e dans
la d\'ecomposition $w=w_1\ldots w_{\kappa(w)}$ comme dans \ref{nu et kappa};
comme la somme des d\'ecalages des orbites de $w_i$  est nulle, il
est produit d'au moins deux cycles $(a_1)_{[h]}$ et $(a_2)_{[k]}$. Posons
$r=(a_1,a_2)$; on a
$(a_1)_{[h]}(a_2)_{[k]}(a_1,a_2)=(a_1,a_2+kn)_{[h+k]}$.
Donc $\nu(wr)=\nu(w)-1$ et $\kappa(wr)=\kappa(w)$, donc $f(wr)=f(w)-1$.

Montrons la propri\'et\'e (2) : Comme $\kappa(w)\leq\nu(w)\leq n$, si $f(w)$ est nul
on doit avoir $\nu(w)=n=\kappa(w)$. La premi\`ere \'egalit\'e prouve que
chaque orbite est un singleton modulo $n$ et la deuxi\`eme prouve alors que
chaque orbite est de d\'ecalage nul, donc $w=1$.

Pour montrer la propri\'et\'e (3) nous montrons que pour $r\in R$ et $w\in W$ on a $f(wr)\leq f(w)+1$. 
Ceci implique par r\'ecurrence que si $w$ est produit de $k$ r\'eflexions on a
$f(w)\leq k$.

Soit $r=(a,b)$ une r\'eflexion. Si $a$ et $b$ apparaissent modulo $n$ dans la m\^eme orbite de
$w$ le lemme \ref{calcul} montre que
que $\nu(wr)=\nu(w)+1$. D'autre part $\kappa(wr)\geq \kappa(w)$, d'o\`u le r\'esultat dans
ce cas. Si $a$ et $b$ sont modulo $n$ dans deux orbites diff\'erentes du m\^eme pseudo-cycle $w_i$
de la d\'ecomposition $w=w_1\ldots w_{\kappa(w)}$, le m\^eme lemme montre
que $\nu(wr)=\nu(w)-1$ et on a
$\kappa(wr)=\kappa(w)$. Si $a$ et $b$ apparaissent dans deux pseudo-cycles distincts,
alors de m\^eme que pr\'ec\'edemment
$\nu(wr)=\nu(w)-1$ et de plus $wr$ a une d\'ecomposition en produit de
$\kappa(w)-1$ \'el\'ements, donc $\kappa(wr)\geq \kappa(w)-1$, ce qui donne bien
$f(wr)\leq
f(w)+1$.
\end{proof}
\begin{corollaire}
Si $w$ est dans un
\quasiparab\  sa longueur de r\'eflexion dans ce sous-groupe est \'egale \`a sa
longueur de r\'eflexion dans $W$.
\end{corollaire}
\begin{proof}On a d\'ej\`a vu que la restriction de $l_R$ \`a un sous-groupe
parabolique est la longueur de r\'eflexion dans ce sous-groupe, gr\^ace au
r\'esultat de Dyer. On peut aussi voir que la formule de la proposition \ref{l_R}
donne la longueur de r\'eflexion dans un parabolique standard. Ceci
donne le r\'esultat dans le cas d'un \quasiparab\  de type $A$.
D'autre part, si un \'el\'ement $w$ est dans le \quasiparab\  de type
$\tilde A_{n-m-1}$ fixant $j+kn$ pour $j$ dans une partie de cardinal $m$
de $[1,n]$ et pour tout $k$, les valeurs de $\nu(w)$
et $\kappa(w)$ diminuent de $m$ dans le \quasiparab\ 
et $n$ est remplac\'e par $n-m$ dans la
formule donnant la longueur. Ceci prouve que $l_R$ se restreint bien aussi
dans ce cas, d'o\`u le r\'esultat.
\end{proof}

Nous allons maintenant chercher quelles r\'eflexions divisent un \'el\'ement de
Coxeter de $W$. La proposition suivante exprime un \'el\'ement de Coxeter comme permutation.
\begin{proposition}\label{coxeter}
Soit $c$ un \'el\'ement de Coxeter de $W$, c'est-\`a-dire le produit des g\'en\'erateurs de Coxeter
dans un ordre arbitraire fix\'e; alors il existe une partition
en deux parties non vides
$$\{1,\ldots,n\}=\{a,b,\ldots,l\}\coprod\{\alpha,\beta,\ldots,\lambda\}$$
telle que $a<b<\ldots<l$ et $\alpha<\beta<\ldots<\lambda$
et que $$c=(a,b,\ldots,l)_{[1]}(\lambda,\ldots,\beta,\alpha)_{[-1]}.$$
\end{proposition}
\begin{proof} L'ensemble des g\'en\'erateurs de Coxeter est
$\{s_1,s_2,\ldots,s_n\}$ avec les notations de \ref{(a,b,c)}.
On a $c=s_{i_1}s_{i_2}\ldots s_{i_n}$, o\`u les indices $i_j$
sont tous distincts.
Quitte \`a changer $c$ en $c\inv$, on peut supposer que dans la suite $(s_{i_1},s_{i_2},\ldots
,s_{i_n})$ l'\'el\'ement $s_1$ est \`a droite de
$s_n$. Alors $s_1$ commute avec tous les $s_i$ qui sont \`a sa droite sauf
\'eventuellement $s_2$ qui \`a son tour commute avec tous les
\'el\'ements qui sont \`a sa droite sauf \'eventuellement $s_3$ \etc\ On peut alors r\'e\'ecrire $c$ comme un produit qui se termine par
$s_1s_2s_3\ldots s_i$
pour un certain $i$. On it\`ere le proc\'ed\'e en commen\c cant avec $s_{i+1}$. Finalement on \'ecrit
$c$ sous la forme
$$(s_{k_h}s_{k_h+1}\ldots s_n)(s_{k_{h-1}}s_{k_{h-1}+1}\ldots
s_{k_h-1})\ldots(s_{k_1}s_{k_1+1}\ldots s_{k_2-1})(s_1s_2\ldots s_{k_1-1}),$$
avec $1< 
k_1<k_2<\ldots<k_h\leq n$
qui est une permutation de la forme voulue: elle s'\'ecrit
$(k_h,k_{h-1},\ldots,k_1)_{[-1]}
(1,2,\ldots,\widehat{k_1}\ldots,\widehat{k_2},\ldots,\widehat{k_h},\ldots,n)_{[1]}$, o\`u $\widehat{k_i}$
signifie que $k_i$ ne figure pas. 
\end{proof}
\begin{notation}\label{X et Xi}
Dans la suite
nous fixons un \'el\'ement de Coxeter $c$ et posons $X=\{a,b,\ldots,l\}+n\bbZ$ et
$\Xi=\{\alpha,\beta,\ldots,\lambda\}+n\bbZ$ comme dans la proposition
pr\'ec\'edente.
\end{notation}

Nous allons d\'emontrer:
\begin{proposition}\label{atomes}
Avec les notations ci-dessus, si $x$ et $y$ sont deux entiers distincts,
la r\'eflexion $(x,y)$ divise $c$ si et seulement
si $x$ et $y$ v\'erifient l'une des propri\'et\'es suivantes:
\begin{enumerate}
\item $x\in X$ et $y\in \Xi$ ou $x\in \Xi$ et $y\in X$.
\item$x, y \in X$, et $|y-x|<n$.
\item$x,y\in \Xi$ et $|y-x|<n$.
\end{enumerate}
\end{proposition}
\begin{proof} On calcule $c(x,y)$ dans tous les cas possibles.

Dans le cas (i), quitte \`a \'echanger les r\^oles de $x$ et $y$ et \`a
changer de repr\'esentants des cycles de $c$,
on peut supposer $x=l$ et $y=\alpha+n$. Le produit $c(x,y)$
vaut $(a,b,\ldots,l,\lambda,\ldots,\beta,\alpha)$ dont la
longueur de r\'eflexion vaut bien $n-1$.

Si $x$ et $y$ sont dans $X$, on peut supposer $x=a$ et on \'ecrit
$c=(a,b,\ldots,t,y+kn,z,\ldots,l)_{[1]}(\lambda,\ldots,\beta,\alpha)_{[-1]}$
pour un certain $k$. Le produit $c(x,y)$ 
vaut $(y,b,\ldots,t)_{[k]}(z,\ldots,l,a+n)_{[1-k]}(\lambda,\ldots,\alpha)_{[-1]}$.
Comme $\nu(c(x,y))=3$, on a
$l_R(c(x,y))=n-1$ si et seulement si
$\kappa(c(x,y))$ est \'egal \`a $2$, c'est \`a dire si on peut regrouper les trois
orbites en deux parties de d\'ecalages nuls  (sinon $\kappa(c(x,y))$ est \'egal \`a 1).
On a donc $\kappa(c(x,y))=2$ si et seulement si $k=0$ ou $k=1$, ce qui donne
le cas (ii).

On fait un raisonnement analogue si $x$ et $y$ sont dans $\Xi$ et on obtient
le cas (iii).
\end{proof}
Les calculs faits dans la preuve pr\'ec\'edente montrent aussi que
\begin{corollaire}\label{quotients} Soit $s$ et $t=n-s$ les cardinaux respectifs
des images de $X$ et $\Xi$ modulo $n$;
un \'el\'ement de longueur $l_R(c)-1$ divise
$c$ si et seulement s'il est d'une des trois formes suivantes:
\begin{enumerate}
\item
$(a_1,a_2,\ldots,a_s,\alpha_t,\alpha_{t-1},\ldots,\alpha_1)$ o\`u
$(a_i)$ et $(\alpha_i)$ sont des sous-suites croissantes form\'ees d'\'el\'ements
cons\'ecutifs respectivement de $X$ et de $\Xi$,
\item
$(a_1,a_2,\ldots,a_r)(a_{r+1},\ldots,a_s)_{[1]}
(\lambda,\ldots,\beta,\alpha)_{[-1]}$
avec $(a_i)$ suite croissante d'\'el\'ements cons\'ecutifs de $X$,
\item
$(a,b,\ldots,l)_{[1]}(\alpha_t,\alpha_{t-1},\ldots,\alpha_{r+1})
(\alpha_r,\alpha_{r-1},\ldots,
\alpha_1)_{[-1]}$ avec $(\alpha_i)$
suite croissante d'\'el\'ements cons\'ecutifs de $\Xi$.
\end{enumerate}
\end{corollaire}

Pour donner la liste des diviseurs de $c$ nous avons besoin de la d\'efinition suivante.
\begin{definition}
On appelle diviseur \'el\'ementaire un pseudo-cycle de $W$ de la forme
$(a_1,\ldots,a_h,\alpha_k,\ldots,\alpha_1)$ avec $h\geq 0$, et $k\geq 0$ et
$h+k\geq 2$, 
ou de la forme
$(a_1,a_2,\ldots,a_h)_{[1]}(\alpha_k,\ldots,\alpha_2,\alpha_1)_{[-1]}$, avec
$h\geq 1$ et $k\geq 1$, o\`u les $a_i$
sont dans $X$ et les $\alpha_i$ dans $\Xi$, et o\`u 
$a_1< a_2< a_3< \ldots< a_h< a_1+n$ et
$\alpha_1<\alpha_2<\ldots<\alpha_k<\alpha_1+n$.
\end{definition}
Remarquons qu'un diviseur \'el\'ementaire est un \'el\'ement de Coxeter d'un
\quasiparab\  de $W$: dans le premier cas il s'agit d'un
\quasiparab\  de type $A_{h+k-1}$, dans le deuxi\`eme d'un
\quasiparab\  de type $\tilde A_{h+k-1}$.

Nous allons donner une interpr\'etation topologique des diviseurs de $c$ en
termes de chemins dans une couronne.
\begin{figure}[htbp]
\begin{center}
 
\begin{picture}(0,0)%
\includegraphics{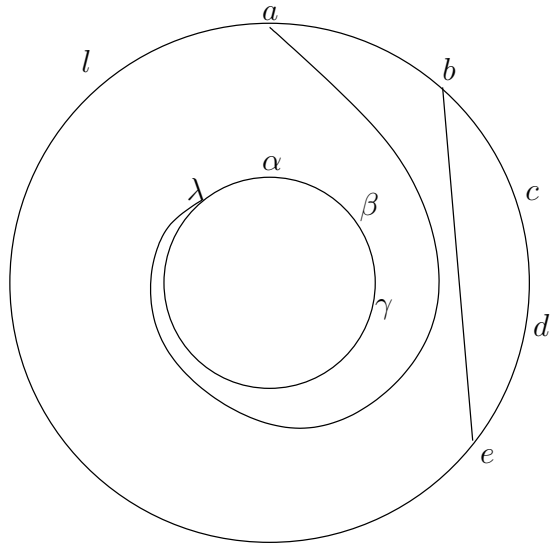}%
\end{picture}%
\setlength{\unitlength}{4144sp}%
\begingroup\makeatletter\ifx\SetFigFont\undefined%
\gdef\SetFigFont#1#2#3#4#5{%
  \reset@font\fontsize{#1}{#2pt}%
  \fontfamily{#3}\fontseries{#4}\fontshape{#5}%
  \selectfont}%
\fi\endgroup%
\begin{picture}(3137,3268)(3839,-6122)
\put(5356,-2986){\makebox(0,0)[lb]{\smash{\SetFigFont{12}{14.4}{\familydefault}{\mddefault}{\updefault}{\color[rgb]{0,0,0}$a$}%
}}}
\put(6436,-3346){\makebox(0,0)[lb]{\smash{\SetFigFont{12}{14.4}{\familydefault}{\mddefault}{\updefault}{\color[rgb]{0,0,0}$b$}%
}}}
\put(6931,-4066){\makebox(0,0)[lb]{\smash{\SetFigFont{12}{14.4}{\familydefault}{\mddefault}{\updefault}{\color[rgb]{0,0,0}$c$}%
}}}
\put(6976,-4876){\makebox(0,0)[lb]{\smash{\SetFigFont{12}{14.4}{\familydefault}{\mddefault}{\updefault}{\color[rgb]{0,0,0}$d$}%
}}}
\put(6661,-5641){\makebox(0,0)[lb]{\smash{\SetFigFont{12}{14.4}{\familydefault}{\mddefault}{\updefault}{\color[rgb]{0,0,0}$e$}%
}}}
\put(5356,-3886){\makebox(0,0)[lb]{\smash{\SetFigFont{12}{14.4}{\familydefault}{\mddefault}{\updefault}{\color[rgb]{0,0,0}$\alpha$}%
}}}
\put(4906,-4066){\makebox(0,0)[lb]{\smash{\SetFigFont{12}{14.4}{\familydefault}{\mddefault}{\updefault}{\color[rgb]{0,0,0}$\lambda$}%
}}}
\put(4276,-3301){\makebox(0,0)[lb]{\smash{\SetFigFont{12}{14.4}{\familydefault}{\mddefault}{\updefault}{\color[rgb]{0,0,0}$l$}%
}}}
\put(5941,-4156){\makebox(0,0)[lb]{\smash{\SetFigFont{12}{14.4}{\familydefault}{\mddefault}{\updefault}{\color[rgb]{0,0,0}$\beta$}%
}}}
\put(6031,-4741){\makebox(0,0)[lb]{\smash{\SetFigFont{12}{14.4}{\familydefault}{\mddefault}{\updefault}{\color[rgb]{0,0,0}$\gamma$}%
}}}
\end{picture}
 
\caption{Deux r\'eflexions sans croisement \label{figure0}}
\end{center}
\end{figure}

On consid\`ere une
couronne dans le plan orient\'e.
On fixe des points \'etiquet\'es $a$, $b$, \dots, $l$
sur le cercle ext\'erieur, dans l'ordre
cyclique et on fait de m\^eme sur le cercle
int\'erieur avec des points \'etiquet\'es
$\alpha$, \dots, $\lambda$.
Cette couronne peut \^etre vue comme le quotient par
la translation de $n$ d'une
bande infinie orient\'ee o\`u les points de $X$ et de $\Xi$ respectivement
sont marqu\'es dans l'ordre croissant
sur chacun des deux bords.
Nous allons associer \`a toute r\'eflexion
une classe de chemins continus dans la couronne.
\begin{notation}\label{representation graphique}
\`A la r\'eflexion $(x,y)$ on associe la classe
d'homotopie \`a extr\'emit\'es fixes de
l'image dans la couronne d'un chemin continu joignant $x$ \`a $y$ dans la
bande (voir figure \ref{figure0},
la repr\'esentation n'est bien d\'efinie qu'une fois fix\'ee
l'identification du quotient de la bande par les translations avec la
couronne).
\end{notation}
Avec cette notation on a:
\begin{corollaire}\label{sans autocroisement}
Une r\'eflexion divise l'\'el\'ement de Coxeter $c$ si et seulement
si elle peut \^etre repr\'esent\'ee par un chemin sans auto-intersection dans la
couronne.
\end{corollaire}
Ce corollaire est une constatation imm\'ediate \`a partir de la liste donn\'ee dans \ref{atomes}.

\begin{figure}[htbp]
\begin{center}
 
\begin{picture}(0,0)%
\includegraphics{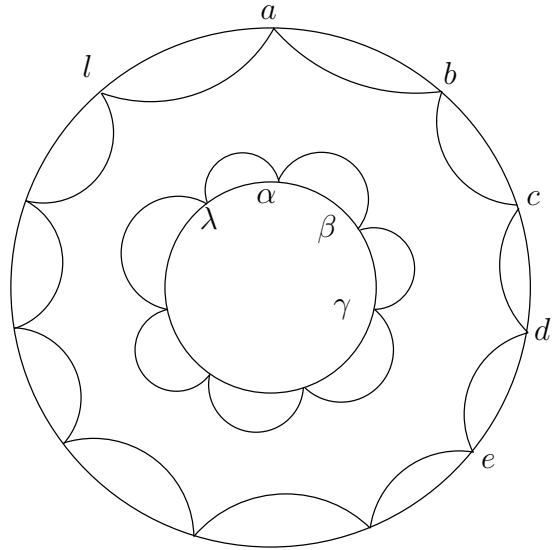}%
\end{picture}%
\setlength{\unitlength}{4144sp}%
\begingroup\makeatletter\ifx\SetFigFont\undefined%
\gdef\SetFigFont#1#2#3#4#5{%
  \reset@font\fontsize{#1}{#2pt}%
  \fontfamily{#3}\fontseries{#4}\fontshape{#5}%
  \selectfont}%
\fi\endgroup%
\begin{picture}(3137,3304)(3839,-6122)
\put(6436,-3346){\makebox(0,0)[lb]{\smash{\SetFigFont{12}{14.4}{\familydefault}{\mddefault}{\updefault}{\color[rgb]{0,0,0}$b$}%
}}}
\put(6931,-4066){\makebox(0,0)[lb]{\smash{\SetFigFont{12}{14.4}{\familydefault}{\mddefault}{\updefault}{\color[rgb]{0,0,0}$c$}%
}}}
\put(6976,-4876){\makebox(0,0)[lb]{\smash{\SetFigFont{12}{14.4}{\familydefault}{\mddefault}{\updefault}{\color[rgb]{0,0,0}$d$}%
}}}
\put(6661,-5641){\makebox(0,0)[lb]{\smash{\SetFigFont{12}{14.4}{\familydefault}{\mddefault}{\updefault}{\color[rgb]{0,0,0}$e$}%
}}}
\put(4276,-3301){\makebox(0,0)[lb]{\smash{\SetFigFont{12}{14.4}{\familydefault}{\mddefault}{\updefault}{\color[rgb]{0,0,0}$l$}%
}}}
\put(5681,-4261){\makebox(0,0)[lb]{\smash{\SetFigFont{12}{14.4}{\familydefault}{\mddefault}{\updefault}{\color[rgb]{0,0,0}$\beta$}%
}}}
\put(5316,-4056){\makebox(0,0)[lb]{\smash{\SetFigFont{12}{14.4}{\familydefault}{\mddefault}{\updefault}{\color[rgb]{0,0,0}$\alpha$}%
}}}
\put(5341,-2951){\makebox(0,0)[lb]{\smash{\SetFigFont{12}{14.4}{\familydefault}{\mddefault}{\updefault}{\color[rgb]{0,0,0}$a$}%
}}}
\put(4976,-4211){\makebox(0,0)[lb]{\smash{\SetFigFont{12}{14.4}{\familydefault}{\mddefault}{\updefault}{\color[rgb]{0,0,0}$\lambda$}%
}}}
\put(5776,-4711){\makebox(0,0)[lb]{\smash{\SetFigFont{12}{14.4}{\familydefault}{\mddefault}{\updefault}{\color[rgb]{0,0,0}$\gamma$}%
}}}
\end{picture}
 
\caption{La repr\'esentation graphique de $\uc$ \label{figure01}}
\end{center}
\end{figure}

Pour pouvoir d\'ecrire tous les diviseurs d'un \'el\'ement de Coxeter nous
utilisons la m\^eme interpr\'etation topologique. \`A
un diviseur \'el\'e\-men\-taire de la forme
$(a_1,a_2,\ldots,a_h,\alpha_k,\ldots,\alpha_2,\alpha_1)$ nous associons 
la classe d'homotopie \`a extr\'emit\'es fixes du lacet compos\'e des
chemins associ\'es aux r\'eflexions
$(a_1,a_2)$, \dots,$(a_{h-1},a_h)$, $(a_h,\alpha_k)$, $(\alpha_k,\alpha_{k-1})$,
\dots, $(\alpha_2,\alpha_1)$, $(\alpha_1,a_1)$ et
\`a un diviseur \'el\'ementaire de la forme
$(a_1,a_2,\ldots,a_h)_{[1]}(\alpha_k,\ldots,\alpha_2,\alpha_1)_{[-1]}$,
nous associons la classe d'homotopie
\`a extr\'emit\'es fixes de l'union
de deux lacets ne se coupant pas et compos\'es respectivement de chemins associ\'es aux r\'eflexions
$(a_1,a_2)$, \dots, $(a_{h-1},a_h)$, $(a_h,a_1+n)$ et $(\alpha_k,\alpha_{k-1})$,
\dots, $(\alpha_2,\alpha_1)$, $(\alpha_1,\alpha_k-n)$ (voir figure \ref{figure01}).

\begin{remarque}\label{sans autointersection}
On peut d\'efinir de la m\^eme fa\c con une interpr\'etation topoloique
de tout pseudo-cycle. Les
diviseurs \'el\'ementaires correspondent exactement
aux pseudo-cycles qui ont une
repr\'esentation (lacet ou union de deux lacets) orient\'ee positivement et
sans autointersection. 
\end{remarque}
\begin{definition}\label{sans croisement}
On dit que deux diviseurs \'el\'ementaires sont sans croisement s'ils admettent des
repr\'esentations sans intersection.
\end{definition}
En particulier si deux diviseurs \'el\'ementaires sont sans croisement, un seul
des deux au plus a des orbites infinies car deux lacets d'image non triviale
dans le groupe fondamental de la couronne et ayant tous deux leur origine sur la
m\^eme composante du bord ont n\'ecessairement une intersection.

Avec cette d\'efinition on peut caract\'eriser les diviseurs d'un \'el\'ement de Coxeter:
\begin{proposition}\label{diviseurs de c}
Les diviseurs de l'\'el\'ement de Coxeter $c$ sont exactement les produits de
diviseurs \'el\'ementaires deux \`a deux sans croisement.
La longueur d'un tel produit est la somme des longueurs des facteurs et ces
facteurs commutent deux \`a deux.
\end{proposition}
La remarque qui pr\'ec\`ede la proposition implique alors
qu'un diviseur de $c$ a z\'ero ou deux orbites infinies (et celles-ci, si elles
existent, sont l'une de
d\'ecalage 1, l'autre de d\'ecalage $-1$).
\begin{proof}
Montrons d'abord l'additivit\'e des longueurs. Dans la formule donnant la
longueur, seul le terme $\kappa$ n'est pas toujours additif. Mais pour un
\'el\'ement qui a 0 ou 2 orbites infinies la d\'ecomposition en pseudo-cycles est unique, un
des pseudo-cycles \'etant form\'e de la permutation induite sur l'union des deux orbites
infinies. Dans ce cas $\kappa$ est bien additif.

Pour prouver l'assertion sur la forme des diviseurs de $c$ nous utilisons le
lemme suivant:

\begin{lemme}\label{reflexions divisant un produit}
Soient $x$ et $y$ deux
entiers distincts et $v$ un \'el\'ement de $W$ ayant 0 ou 2
orbites infinies.
\begin{itemize}
\item[(i)]
La r\'eflexion $(x,y)$ divise $v$ si et seulement si elle divise un des
pseudo-cycles
de $v$.
\item[(ii)] La r\'eflexion $(x,y)$ divise un pseudo-cycle si et seulement si $x$
et $y$ sont dans la m\^eme orbite de ce pseudo-cycle ou si chacun d'eux
est dans une des orbites infinies de ce pseudo-cycle.
\end{itemize}
\end{lemme}
\begin{proof}
Si $x$ est dans une orbite finie de $v$ et $y$ dans
une orbite finie ou infinie distincte de celle de $x$,
alors le calcul de \ref{calcul}, appliqu\'e avec $a_1=y$, $a_i=x$ et $k=0$, et en faisant passer $(x,y)$
dans le membre de droite, montre que
$\nu(v(x,y))=\nu(v)-1$  et que la suite des d\'ecalages non nuls des orbites est la m\^eme pour
$v$ et $v(x,y)$.
Donc $\kappa(v(x,y))=\kappa(v)-1$ et
$l_R(v(x,y))=l_R(v)+1$, donc $(x,y)$ ne divise pas $v$.
R\'eciproquement, si $x$ et $y$ sont dans le support
d'un m\^eme pseudo-cycle de $v$, \ref{calcul} appliqu\'e avec $a_1=x$, $a_2=y$ et $h=0$ et avec ou bien $k=1$
ou bien $k=0$ montre que $(x,y)$ divise ce pseudo-cycle.
Comme les longueurs des pseudo-cycles de $v$
s'ajoutent d'apr\`es le d\'ebut de la d\'emonstration de la proposition, 
on en d\'eduit que $(x,y)$ divise $v$.
\end{proof}

Nous prouvons maintenant, par
r\'ecurrence descendante sur la longueur du diviseur,
l'assertion de la proposition \ref{diviseurs de c} sur les diviseurs de $c$. Il y a
exactement un diviseur de longueur $n$ qui est $c$ lui-m\^eme.
Un diviseur de longueur $k$ est le produit d'un diviseur $v$ de longueur $k+1$ par
une r\'eflexion $(x,y)$ qui divise $v$. Par hypoth\`ese de r\'ecurrence
appliqu\'ee \`a $v$, les pseudo-cycles de $v$ sont des diviseurs \'el\'ementaires deux \`a deux
sans croisement.
D'apr\`es le lemme pr\'ec\'edent la r\'eflexion $(x,y)$ divise
$v$ si et seulement
si elle divise un des pseudo-cycles $v_1$ de $v$.
On est ramen\'e \`a trouver tous les diviseurs de longueur $l_R(v_1)-1$
d'un diviseur \'el\'ementaire $v_1$.
Si $v_1$ a deux orbites infinies,
c'est un \'el\'ement analogue \`a $c$ dans le sous-groupe
\quasiparab\  correspondant au support de $v_1$. Les diviseurs cherch\'es
sont donn\'es par \ref{quotients}. Ce sont bien des produits de diviseurs
\'el\'ementaires deux \`a deux sans croisement, de supports inclus dans le support de
$v_1$, donc aussi sans croisement avec les autres pseudo-cycles de
$v$.
Si $v_1$ n'a pas d'orbite infinie il s'\'ecrit
$(a_1,\ldots,a_h)$  o\`u $l_R(v_1)=h-1$ et ses diviseurs de longueur $h-2$
sont exactement les \'el\'ements
$(a_1,a_2,a_{i-1},a_{j+1},\ldots,a_h)(a_i,a_{i+1},a_{i+2},\ldots,a_j)$,
ce qui est aussi de la forme annonc\'ee.

R\'eciproquement, si $v$ est un produit de diviseurs \'el\'ementaires deux \`a deux
sans croisement et si $l_R(v)<n$ on va montrer qu'on peut trouver une
r\'eflexion $r$ telle que $r v$ soit de longueur $l_R(v)+1$
et soit un produit de diviseurs \'el\'ementaires deux \`a deux sans croisement.
L'hypoth\`ese de r\'ecurrence montre alors que $r v$ divise $c$, donc que $v$
divise $c$. Distinguons plusieurs cas pour $v$.
\begin{enumerate}
\item[$\bullet$]
Si $v$ a une orbite infinie, il en a alors exactement deux, d'apr\`es la
remarque qui suit \ref{sans croisement} et ces orbites sont de d\'ecalages $1$
et $-1$, donc les cycles correspondants sont de la forme
$(a_1,a_2,\ldots,a_h)_{[1]}$ et
$(\alpha_k,\ldots,\alpha_2,\alpha_1)_{[-1]}$ avec
$a_1< a_2< a_3< \ldots< a_h< a_1+n$ et
$\alpha_1<\alpha_2<\ldots<\alpha_k<\alpha_1+n$, les $a_i$ \'etant dans $X$ et
les $\alpha_i$ dans $\Xi$;
d'autre part $v$ doit avoir aussi une orbite finie car $v\neq c$. Tous les points
d'une telle orbite finie sont compris entre deux \'el\'ements cons\'ecutifs 
d'une des deux orbites infinies, puisque les diviseurs \'el\'ementaires sont sans croisement.
On peut supposer que cette orbite est comprise dans $[a_1,a_2]$. Notons-la
$\{b_1,\ldots,b_l\}$ avec $a_1<b_1<b_2<\cdots<b_l<a_2$. On a alors 
\begin{multline*}
(a_2,b_1)
(\alpha_k,\ldots,\alpha_2,\alpha_1)_{[-1]}
(a_1,a_2,\ldots,a_h)_{[1]}(b_1,\ldots,b_l)=\\
(\alpha_k,\ldots,\alpha_2,\alpha_1)_{[-1]}
(a_1,b_1,\ldots,b_l,a_2,\ldots,a_h)_{[1]}.
\end{multline*}
\end{enumerate}

 Si $v$ n'a pas d'orbite infinie, il a au moins une orbite
finie de la forme $(a_1,a_2\ldots,a_h,\alpha_k,\ldots,\alpha_2,\alpha_1)$ avec les m\^emes
conventions que pr\'ec\'edemment (cette orbite peut \^etre un singleton).
Quitte \`a \'echanger les r\^oles de $X$ et de $\Xi$
on peut supposer que $h\neq 0$. Soit $b_1$ le successeur de $a_h$ dans $X$.
Il y a trois cas:
\begin{enumerate}
\item[$\bullet$]
Si $b_1=a_1+n$ alors 
\begin{multline*}
(b_1,\alpha_k)
(a_1,a_2\ldots,a_h,\alpha_k,\ldots,\alpha_2,\alpha_1)=\\
(a_1,a_2,\ldots,a_h)_{[1]}
(\alpha_k,\ldots,\alpha_2,\alpha_1)_{[-1]}.
\end{multline*}
\item[$\bullet$]
Si $k\neq 0$ et si le cycle de $v$ dont le support contient $b_1$ est de la forme
$(b_1,b_2\ldots,b_l,\beta_m,\ldots,\beta_2,\beta_1)$
avec $b_1< b_2<\ldots<b_l< b_1+n$ et
$\beta_1<\beta_2<\ldots<\beta_m<\beta_1+n$, les $b_i$ \'etant dans $X$ et
les $\beta_i$ dans $\Xi$; alors
\begin{multline*}
(b_1,\alpha_k)(a_1,a_2,\ldots,a_h,\alpha_k,\ldots,\alpha_2,\alpha_1)
(b_1,b_2,\ldots,b_l,\beta_m,\ldots,\beta_2,\beta_1)=\\
(a_1,a_2,\ldots,a_h,b_1,b_2,\ldots,b_l,\beta_m,\ldots,\beta_2,\beta_1,
\alpha_k,\ldots,\alpha_2,\alpha_1).
\end{multline*}
\item[$\bullet$]
Si $k=0$ et si le cycle de $v$ dont le support contient de $b_1$ est de la forme
\[(b_1,b_2\ldots,b_r,\beta_m,\ldots,\beta_2,\beta_1,b_{r+1},\ldots,b_l)\]
avec $b_{r+1}<b_{r+2}<\cdots<b_l<a_1<a_h<b_1< b_2<\ldots<b_r< b_{r+1}+n$ et
$\beta_1<\beta_2<\ldots<\beta_m<\beta_1+n$, les $b_i$ \'etant
dans $X$ et les $\beta_i$ dans $\Xi$; alors
\begin{multline*}
(a_1,b_1)(a_1,a_2,\ldots,a_h)
(b_1,b_2\ldots,b_r,\beta_m,\ldots,\beta_2,\beta_1,b_{r+1},\ldots,b_l)=\\
(a_1,a_2,\ldots,a_h,b_1,b_2,\ldots,b_r,\beta_m,\ldots,\beta_2,\beta_1,
b_{r+1},\ldots,b_l).
\end{multline*}
\end{enumerate}
Dans tous les cas le produit est bien comme annonc\'e.

Enfin puisque les longueurs de diviseurs \'el\'ementaires deux \`a deux
sans croisements s'ajoutent et que de tels diviseurs commutent dans $W$, ils commutent aussi
dans le mono{\"\i}de.
\end{proof}
\begin{remarque}\label{non-additivite}
Il est faux que la longueur d'un produit de pseudo-cycles de supports disjoints est
la somme des longueurs des pseudo-cycles, comme le montre l'exemple
suivant dans $\tilde A_n$ avec $n\geq 5$:
\[w=[(1)_{[-1]}(2)_{[-1]}(3)_{[2]}][(4)_{[1]}(5)_{[1]}(6)_{[-2]}].\]
La longueur de chacun des facteurs
est \'egale \`a 4 ($\nu=3$, $\kappa=1$, dans un \quasiparab\  de type $\tilde
A_2$) et la longueur du produit
vaut 6 car le m\^eme \'el\'ement s'\'ecrit
\[w=[(1)_{[-1]}(4)_{[1]}][(2)_{[-1]}(5)_{[1]}][(3)_{[2]}(6)_{[-2]}],\]
ce qui prouve que $\kappa$ vaut $3$ (et on a $\nu=6$) dans un \quasiparab\ de
type $\tilde A_5$.

Par contre s'il y a au plus deux orbites infinies on a bien
additivit\'e des longueurs.
\end{remarque}
\`A tout diviseur $w$ de $c$ on peut associer la partition p\'eriodique
de $\bbZ$ dont les parties sont les supports des pseudo-cycles
de $w$.

Nous dirons que deux parties $A$ et $B$ de $\bbZ$
sont sans croisement
si pour tous $x$ et $y$ dans $A$ et tous $z$ et $t$ dans $B$ il existe deux
chemins sans intersection dans la bande joignant respectivement $x$ \`a $y$ et
$z$ \`a $t$.
On peut alors r\'eexprimer \ref{diviseurs de c} et \ref{sans autointersection} par:
\begin{corollaire} \label{partition}
Une partition p\'eriodique de $\bbZ$ 
dont toute partie infinie rencontre \`a la fois $X$ et $\Xi$
est associ\'ee comme ci-dessus \`a un diviseur (unique) de $c$ si et seulement si
ses parties sont deux \`a
deux sans croisement.
\end{corollaire}
\section{Les mono{\"\i}des duaux}\label{monoide dual}
L'objectif de cette partie est de prouver la conjecture \ref{conjecture} (ii)
pour le type $\tilde A$.
Nous suivons une d\'emarche analogue \`a celle de \cite{B} (ou de \cite{BDM}).
En particulier nous nous pla\c cons dans le cadre des groupes positivement
engendr\'es telle qu'elle est
expos\'ee dans \cite[0.4 et 0.5]{B}. 
Rappelons-en les r\'esultats principaux
dans un cadre plus g\'en\'eral car nous ne
supposons pas que le nombre de g\'en\'erateurs est fini.
Soit $G$ un groupe engendr\'e comme mono{\"\i}de
par un ensemble $R$. On dit que $(G,R)$ est un groupe
positivement engendr\'e. On d\'efinit
la longueur dans $l_R(g)$ par rapport \`a $R$ de $g\in G$ comme le
nombre minimum de facteurs dans une d\'ecomposition de $g$ en produit de
g\'en\'erateurs. On dit que $h\in G$ divise $g\in G$ \`a gauche  et que $k\in G$
divise $g$ \`a droite si $g=hk$ avec $l_R(g)=l_R(h)+l_R(k)$.
On dit qu'un \'el\'ement est \'equilibr\'e si ses diviseurs \`a droite et \`a gauche sont
les m\^emes.
\begin{definition}
Soit $(G,R)$ un groupe positivement
engendr\'e, soit $c\in G$ un \'el\'ement \'equilibr\'e et soit
$P_c$ l'ensemble des diviseurs (\`a
gauche ou \`a droite) de $c$. On consid\`ere un ensemble
$\uP_c=\{\,\underline w\,\mid\,w\in P_c\}$
en bijection avec $P_c$ et on d\'efinit un mono{\"\i}de not\'e
$M(P_c)$ par la pr\'esentation suivante: l'ensemble des g\'en\'erateurs est $\uP_c$ 
et les relations sont $\underline w.\underline w'=\underline{ww'}$ pour tous les couples
$(w,w')$ tels que $w$, $w'$ et $ww'$
sont des diviseurs de $c$ et que $l_R(ww')=l_R(w)+l_R(w')$.
\end{definition}
Les notions standard de divisibilit\'e \`a gauche ou \`a droite
dans le mono\"\i de $M(P_c)$ \'etendent les notions correspondantes
de divisibilit\'e d\'efinies dans
$P_c$. Notons que les atomes, \ie,
les \'el\'ements diff\'erents de 1 qui ne sont pas produit de deux facteurs
diff\'erents de 1 sont des \'el\'ements de $P_c$.

Remarquons qu'on a un morphisme de mono{\"\i}des $M(P_c)\rightarrow G$ donn\'e par $\uw\mapsto
w$ pour $w\in P_c$. L'existence de ce morphisme permet facilement
de prouver que $M(P_c)$ a une propri\'et\'e de simplifiabilit\'e partielle \`a gauche
et \`a droite (\cf\ \cite[0.4.4]{B}): si $am=bm$ ou si $ma=mb$ avec $a$ et $b$
dans $P_c$ et $m\in M(P_c)$ alors $a=b$. On a de plus:
\begin{proposition}\label{conjugaison par c}
Pour tout $w\in P_c$ l'\'el\'ement $w'=cwc\inv$ est dans $P_c$ et est l'unique
\'el\'ement de $P_c$
tel que $\uw'.\uc=\uc.\uw$. L'application $\uw\mapsto\uw'$ d\'efinit un automorphisme
du mono{\"\i}de $M(P_c)$.
\end{proposition}
On  appellera ``conjugaison par $\uc$'' l'automorphisme ainsi d\'efini.
\begin{proof}
On a $c=xw$ avec $l_R(c)=l_R(w)+l_R(x)$.
Comme $c$ est \'equilibr\'e, l'\'el\'ement
$x$ est aussi un diviseur de $c$ \`a droite, donc on peut \'ecrire $c=w'x$ avec
$l_R(w')+l_R(x)=l_R(c)$. On a $w'=cwc\inv$ et
$\uc=\uw'.\ux=\ux.\uw$, d'o\`u $\uc.\uw=\uw'.\ux.\uw=\uw'.\uc$.
L'application $\uw\mapsto\uw'$
d\'efinit un morphisme de mono{\"\i}des. On a un morphisme en sens
inverse en faisant un raisonnement analogue en partant de $w'$. D'o\`u le
r\'esultat.
\end{proof}
Notons encore la propri\'et\'e g\'en\'erale suivante de $M(P_c)$:
\begin{proposition}
Tout \'el\'ement de $M(P_c)$  divise une puissance suffisament grande de $c$.
\end{proposition}
\begin{proof}
Tout \'el\'ement de $M(P_c)$ s'\'ecrit
$\uw_1\ldots \uw_k$ pour un certain $k$, o\`u les $w_i$ sont des \'el\'ements de
$P_c$.
On montre par r\'ecurrence sur $k$ que
$\uw_1\ldots \uw_k$ divise $\uc^k$: par d\'efinition les \'el\'ements de $\uP_c$
divisent $\uc$, donc il existe $x\in P_c$ tel que
$\uw_k\ux=\uc$. Par \ref{conjugaison par c} on a
$\uw_1\ldots\uw_k\ux=\uc\uw'_1\ldots\uw'_{k-1}$ o\`u les $w'_i$ sont des
\'el\'ements de $P_c$. Par
hypoth\`ese de r\'ecurrence $\uw'_1\ldots\uw'_{k-1}$ divise $\uc^{k-1}$, d'o\`u la
proposition.
\end{proof}
Remarquons que par d\'efinition
la longueur $l_R$ s'\'etend en une fonction additive sur $M(P_c)$.
Remarquons aussi que les diviseurs de $\uc$ dans $M(P_c)$ sont exactement les
\'el\'ements de $\uP_c$. Les \'el\'ements $\ur\in\uP_c$ tels que $r\in R$ sont les
atomes du mono{\"\i}de $M(P_c)$.

Revenons \`a la situation du groupe de Coxeter $W$ de type $\tilde A_{n-1}$ et
appliquons les constructions pr\'ec\'edentes \`a
un \'el\'ement de Coxeter $c$ fix\'e de
$W$ comme dans la section pr\'ec\'edente, dont
on garde les notations. 

Nous prouvons maintenant dans ce cas la conjecture \ref{conjecture} (ii):
\begin{proposition}\label{presentation}
Le mono{\"\i}de $M(P_c)$ est engendr\'e par les $\ur$ o\`u $r$ est une r\'eflexion qui divise $c$
avec comme relations
\begin{equation}
\ur.\ut=\underline{rtr}.\ur\label{*}\end{equation}
si $r$ et $t$ sont deux r\'eflexions distinctes telles que $rt$ divise $c$.
\end{proposition}
Remarquons qu'un cas particulier de ces relations est que $\ur.\ut=\ut.\ur$
si $rt$ divise $c$ et que $r$ et $t$ commutent.
\begin{proof}
La preuve suit les m\^emes grandes lignes que celle de \cite[2.1.4]{B}.
Le mono{\"\i}de $M(P_c)$ est engendr\'e par les $\ur$ o\`u $r$ est une r\'eflexion qui
divise $c$ et les relations \ref{*} sont vraies dans $M(P_c)$. Il suffit
de voir que ces relations impliquent les autres, c'est-\`a-dire que pour tout
$w\in P_c$ on peut passer d'une \'ecriture de $w$ de longueur $l_R(w)$ \`a une
autre uniquement par les relations
\begin{equation}
r.s=(rsr).r\label{**}\end{equation}
si $r$ et $s$ sont deux r\'eflexions distinctes telles que $rs$ divise $c$.
Prouvons ceci par r\'ecurrence sur $l_R(w)$.
Si $l_R(w)=1$ il n'y a qu'une \'ecriture de longueur minimale de $w$. Dans le
cas g\'en\'eral
il suffit de prouver que si $t\in R$ et $t\preccurlyeq w\in P_c$
alors \`a partir d'une \'ecriture minimale de $w$ fix\'ee et
par application des relations \ref{**}
on peut obtenir une \'ecriture minimale de $w$
commen\c cant par $t$. L'hypoth\`ese de r\'ecurrence permet alors de
conclure.
L'\'el\'ement $w$, divisant $c$, est un produit de diviseurs
\'el\'ementaires comme dans \ref{diviseurs de c}.
On fixe une \'ecriture minimale de $w$ obtenue par concat\'enation
d'une \'ecriture minimale de chaque diviseur \'el\'ementaire.
Par \ref{reflexions divisant un produit}
une r\'eflexion $t$ divise $w$ si et seulement si elle divise
un des diviseurs \'el\'ementaires de
$w$.
Comme deux r\'eflexions qui interviennent dans les \'ecritures de deux diviseurs
\'el\'ementaires de $w$ distincts commutent entre elles,
$t$ commute avec les r\'eflexions qui interviennent dans les \'ecritures
minimales des autres diviseurs \'el\'ementaires. On est donc ramen\'e \`a montrer le
r\'esultat pour un seul diviseur \'el\'ementaire. Si ce diviseur est un \'el\'ement de
Coxeter d'un groupe de type $A$ (c'est-\`a-dire a toutes ses orbites finies) le
r\'esultat est connu (\cf\ \cite{B} et \cite{BDM}). On est ramen\'e au cas d'un
\'el\'ement de Coxeter d'un groupe de type $\tilde A$. Il suffit donc de prouver
le r\'esultat pour $c$ lui-m\^eme. Quitte \`a faire une permutation circulaire
des r\'eflexions \'el\'ementaires, ce qui revient \`a une conjugaison donc laisse
invariantes les relations \ref{**}, on peut ramener $c$ \`a \^etre de la forme
(\cf\ preuve de \ref{coxeter})
\[c=(s_{k_h}s_{k_h+1}\ldots s_{k_{h+1}-1})(s_{k_{h-1}}s_{k_{h-1}+1}\ldots
s_{k_h-1})\ldots(s_1s_2\ldots s_{k_1-1}),\]
avec $1=k_0<k_1<k_2<\ldots<k_h<k_{h+1}=n+1$ et $s_i=(i,i+1)$.
On part de cette \'ecriture et on veut faire
appara{\^\i}tre $t$ \`a gauche de $c$ par application des relations \ref{**}.
En fait il suffit de faire appara{\^\i}tre $t$ dans une \'ecriture de $c$; on peut
ensuite le ramener \`a gauche par application des relations \ref{**}. 
Chaque $(s_{k_i}s_{k_i+1}\ldots s_{k_{i+1}-1})$ est un \'el\'ement de Coxeter d'un
groupe de type $A$, donc on peut faire appara{\^\i}tre \`a gauche ou \`a droite, par
application des relations \ref{**} n'importe quelle r\'eflexion de support inclus
dans $[k_i,k_{i+1}]$. On en d\'eduit que si $k_{i-1}\leq a\leq k_i\leq k_j\leq
b\leq k_{j+1}$, on peut faire appara{\^\i}tre
dans l'\'ecriture de $c$  par application de \ref{**} le produit
$(k_j,b)(k_{j-1},k_j)\ldots(k_{i+1},k_i)(k_i,a)$ qui est l'\'ecriture d'un
\'el\'ement de Coxeter d'un groupe de type $A$, donc on peut faire appara{\^\i}tre
$(a,b)$ dans l'\'ecriture de cet \'el\'ement. Le m\^eme type d'argument 
montre qu'on peut faire appara{\^\i}tre $(b,a+n)$ sous les m\^emes hypoth\`eses,
en faisant appara{\^\i}tre le produit $(b,1+n)$ \`a droite du produit
$(s_{k_h}\ldots s_{k_{h+1}})\ldots(s_{k_j}\ldots s_{k_{j+1}})$
et $(1+n,a+n)=(1,a)$ \`a gauche du produit 
$(s_{k_i}\ldots s_{k_{i+1}-1})\ldots(s_{k_0}\ldots s_{k_1-1})$.

Il reste \`a voir qu'on peut faire appara{\^\i}tre toute r\'eflexion $t=(a,\alpha)$
o\`u $a\in X$ et $\alpha\in \Xi$. Comme
$c$ est la translation d'une position dans le sens croissant de $X$ et d'une position dans
le sens d\'ecroissant de $\Xi$, il
conjugue $(a,\alpha)$ sur $(a',\alpha')$ o\`u $a'$ est translat\'e de $a$ dans $X$
d'une position dans le sens croissant et $\alpha'$ est translat\'e de $\alpha$ dans $\Xi$
d'une position dans le sens d\'ecroissant. Donc $(a,\alpha)$ peut \^etre
ramen\'e par conjugaison par une puissance de $c$ sur une r\'eflexion d'une des formes
$(a,b)$ ou $(a,b-n)$, avec $1\leq a,b\leq n$, et par la premi\`ere partie de cette
d\'emonstration on sait qu'on peut faire appara{\^\i}tre une telle r\'eflexion
dans une \'ecriture de $c$ par application des relations \ref{**}.
La conjugaison par $c$ est aussi une suite
d'applications de \ref{**}, d'o\`u le r\'esultat.
\end{proof}

\section{Pr\'esentations duales pour les groupes d'Artin-Tits affines de type $\tilde A$}
Nous gardons les notations des deux sections pr\'ec\'edentes, en particulier
$c$ est un \'el\'ement de Coxeter du groupe de Coxeter $W$ de type
$\tilde A_{n-1}$.
Le premier but de cette section est de prouver le th\'eor\`eme suivant (\cf\ 
conjecture \ref{conjecture} (i)):
\begin{theoreme}\label{G isom B}
Le groupe de fractions $G(P_c)$ de $M(P_c)$
est isomorphe au groupe d'Artin-Tits de type $\tilde
A_{n-1}$.
\end{theoreme}
Nous noterons $\Bn$ le groupe d'Artin-Tits de type $\tilde A_{n-1}$.
Il est engendr\'e par $\bs_1,\bs_2,\ldots \bs_n$, avec comme relations
$\bs_i\bs_{i+1}\bs_i=\bs_{i+1}\bs_i\bs_{i+1}$ pour
$i=1,\ldots,n$, si on pose $\bs_{n+1}=\bs_1$, et $\bs_i\bs_j=\bs_j\bs_i$ si
$i\neq j\pm 1\pmod n$.
Pour prouver le th\'eor\`eme \ref{G isom B} on montre d'abord que dans $G(P_c)$ les
\'el\'ements $\us_1,\ldots,\us_n$ v\'erifient les m\^emes relations de tresses.
Ceci d\'efinit un morphisme $\Bn\to G(P_c)$. On trouve ensuite dans
$\Bn$ des \'el\'ements dont les images par ce morphisme sont les g\'en\'erateurs $\ur$ de
$G(P_c)$ et qui v\'erifient les relations \ref{*}, ce qui prouve
la bijectivit\'e du morphisme.
Pour cette deuxi\`eme \'etape on utilisera l'interpr\'etation de $\Bn$
comme groupe fondamental.

\begin{proposition}\label{B to G}
L'application $\bs_i\mapsto\us_i $  pour $i=1,\ldots,n$
se prolonge en un homomorphisme 
$\Bn\rightarrow G(P_c)$.
\end{proposition}
\begin{proof}
Cela revient \`a montrer que  dans $G(P_c)$ on a 
$\us_i\us_{i+1}\us_i=\us_{i+1}\us_i\us_{i+1}$
pour $i=1,\ldots,n$ et
$\us_i\us_j=\us_j\us_i$ si $i,j\in[1,n]$ et $|i-j|\geq 2$.
Par \ref{presentation} on a
$\us_i\us_{i+1}=\us'\us_i=\us_{i+1}\us'$,
o\`u $s'=(i,i+2)$.
On en d\'eduit
$\us_{i+1}\us_i\us_{i+1}=\us_{i+1}\us'\us_i=\us_i\us_{i+1}\us_i$.
On a aussi par \ref{presentation} 
$\us_i\us_j=\us_j\us_i$ si $|i-j|\geq 2$. 
\end{proof}

Avant de montrer que ce morphisme est un
isomorphisme, nous rappelons l'interpr\'etation de $\Bn$ comme
sous-groupe du groupe de tresses
\`a $n$ brins dans $\bbC^*$ (\cf\  \cite{Graham-Lehrer} et
\cite{Allcock}).
On consid\`ere un $n$-uplet de points de $\bbC^*$.
Le groupe des tresses dans $\bbC^*$ de base ce $n$-uplet (``tresses \`a $n$
brins'' dans $\bbC^*$), est isomorphe au
groupe d'Artin-Tits $B(B_n)$ de type $B_n$. L'application qui associe \`a une
telle tresse le nombre de tours total des brins autour de 0 est un morphisme
\`a valeur dans $\bbZ$. 
Le groupe des tresses de type $\tilde A_{n-1}$
est le noyau de ce morphisme. Notons aussi qu'on peut consid\'erer 
le groupe des tresses \`a $n$
brins dans $\bbC^*$ de base le $n$-uplet $(x_1,\ldots,x_n)$ comme le
sous-groupe du groupe des tresses \`a $n+1$ brins dans $\bbC$ de base le
$n+1$-uplet $(0,x_1,\ldots,x_n)$ tel que le brin issu de l'origine soit
trivial (tresses pures relativement \`a un brin fix\'e). 

Nous nous pla\c cons dans le cadre de
\ref{representation graphique} et
nous choisissons comme $n$-uplet
de base $(a,b,\ldots,l,\alpha,\ldots,\lambda)$ comme dans
\ref{representation graphique}.
Nous associons \`a chaque
$\underline r$ o\`u $r$ est une r\'eflexion de
$P_c$ une tresse de $\Bn$ la fa\c con suivante: $r$ est repr\'esent\'ee
par un chemin $\gamma$ dans la couronne, sans auto-intersection, reliant 
$i$ \`a $j$ o\`u $i$  et $j$ sont les images dans la couronne de deux
points de $X\cup\Xi$.
On associe \`a ce chemin la tresse o\`u tous les points sont
fixes sauf les deux points partant respectivement de $i$ et $j$ qui
suivent $\gamma$ en sens inverse et s'``\'evitent par la droite'' si
l'orientation du plan est choisie dans le sens horaire ce que nous supposerons
dans les figures qui suivent.
Plus pr\'ecis\'ement, on peut supposer que $\gamma$ est une application diff\'erentiable
de $[0,1]$ dans la couronne telle que les tangentes en 0 et 1 soient orthogonales au bord de la couronne. Soit $\vec n(t)$ un vecteur normal dans le sens
direct \`a $\gamma$ en $\gamma(t)$. On consid\`ere la tresse o\`u tous les brins
sont fixes sauf un brin partant de $i$ donn\'e par
$t\mapsto\gamma(t)+\varepsilon\sin(\pi t)\vec n(t)$ et un
brin partant de $j$ donn\'e
par $t\mapsto\gamma(1-t)-\varepsilon\sin(\pi (1-t))\vec n(1-t)$ o\`u
$\varepsilon$ est assez petit pour que la tresse soit dans la couronne.

\begin{proposition}\label{G to B}
L'application que nous venons de d\'efinir se prolonge en un isomorphisme
de $G(P_c)$ dans $\Bn$
inverse de l'homomorphisme d\'efini par \ref{B to G}.
\end{proposition}
\begin{proof}
Il faut voir que les relations \ref{*}
sont v\'erifi\'ees par les images des
\'el\'ements $\ur$. On d\'eduit de \ref{diviseurs de c} qu'il y a
trois types de couples $(r,t)\in R^2$ tels que $rt$ divise $c$.
\begin{enumerate}
\item[$\bullet$]
Si $r$ et $t$ correspondent \`a des chemins sans intersection $\ur$ et $\ut$
commutent et il est clair que leurs images dans $\Bn$ commutent
aussi.
\item[$\bullet$]
Si $r=(i,j)$ et $t=(j,k)$ o\`u $i$, $j$ et $k$ sont deux \`a deux distincts modulo
$n$ et si $i,j,k$ sont les sommets d'un triangle curviligne
direct dans la bande, la relation est
$\underline{(i,j)}.\underline{(j,k)}=\underline{(i,k)}.\underline{(i,j)}$.
Cette relation est v\'erifi\'ee par les tresses images:
si nous notons encore $i$, $j$ et $k$ les images repectives
dans la couronne des points $i$, $j$ et $k$, les \'el\'ements
$\underline{(i,j)}$, $\underline{(j,k)}$ et
$\underline{(i,k)}$
correspondent \`a des chemins respectivement de $i$ \`a $j$
de $j$ \`a $k$ et de $k$ \`a $i$ formant le bord d'un triangle curviligne direct inclus
dans la couronnne et la
relation pour les tresses correspondantes n'est autre que
la relation classique pour les
tresses \`a trois brins (voir figure \ref{relation standard}) dans le groupe
de tresses de $\bbC^*$ (ou de $\bbC$).

\begin{figure}[htbp]
\begin{center}
 
\begin{picture}(0,0)%
\includegraphics{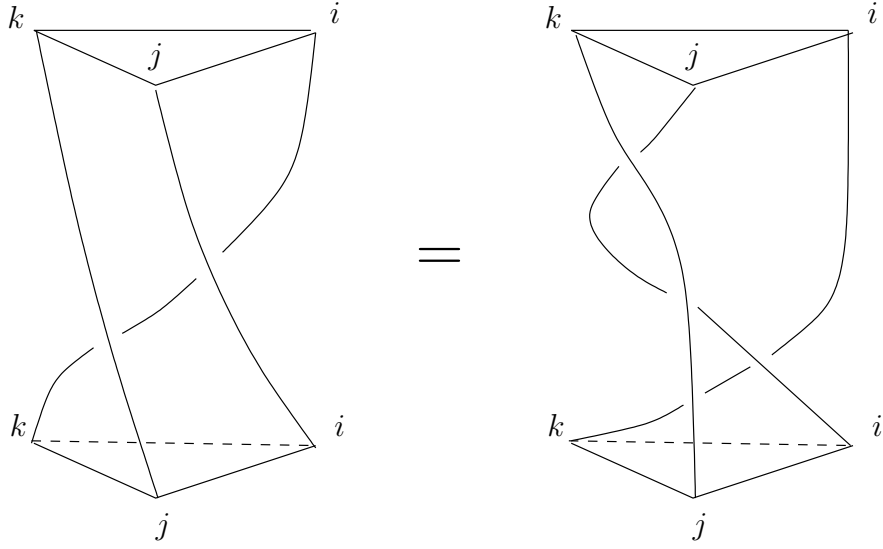}%
\end{picture}%
\setlength{\unitlength}{3947sp}%
\begingroup\makeatletter\ifx\SetFigFont\undefined%
\gdef\SetFigFont#1#2#3#4#5{%
  \reset@font\fontsize{#1}{#2pt}%
  \fontfamily{#3}\fontseries{#4}\fontshape{#5}%
  \selectfont}%
\fi\endgroup%
\begin{picture}(5430,3399)(1771,-7240)
\put(3796,-3976){\makebox(0,0)[lb]{\smash{\SetFigFont{12}{14.4}{\familydefault}{\mddefault}{\updefault}{\color[rgb]{0,0,0}$i$}%
}}}
\put(3826,-6571){\makebox(0,0)[lb]{\smash{\SetFigFont{12}{14.4}{\familydefault}{\mddefault}{\updefault}{\color[rgb]{0,0,0}$i$}%
}}}
\put(2701,-7186){\makebox(0,0)[lb]{\smash{\SetFigFont{12}{14.4}{\familydefault}{\mddefault}{\updefault}{\color[rgb]{0,0,0}$j$}%
}}}
\put(2656,-4231){\makebox(0,0)[lb]{\smash{\SetFigFont{12}{14.4}{\familydefault}{\mddefault}{\updefault}{\color[rgb]{0,0,0}$j$}%
}}}
\put(1771,-4006){\makebox(0,0)[lb]{\smash{\SetFigFont{12}{14.4}{\familydefault}{\mddefault}{\updefault}{\color[rgb]{0,0,0}$k$}%
}}}
\put(1786,-6571){\makebox(0,0)[lb]{\smash{\SetFigFont{12}{14.4}{\familydefault}{\mddefault}{\updefault}{\color[rgb]{0,0,0}$k$}%
}}}
\put(7171,-3976){\makebox(0,0)[lb]{\smash{\SetFigFont{12}{14.4}{\familydefault}{\mddefault}{\updefault}{\color[rgb]{0,0,0}$i$}%
}}}
\put(7201,-6571){\makebox(0,0)[lb]{\smash{\SetFigFont{12}{14.4}{\familydefault}{\mddefault}{\updefault}{\color[rgb]{0,0,0}$i$}%
}}}
\put(6076,-7186){\makebox(0,0)[lb]{\smash{\SetFigFont{12}{14.4}{\familydefault}{\mddefault}{\updefault}{\color[rgb]{0,0,0}$j$}%
}}}
\put(6031,-4231){\makebox(0,0)[lb]{\smash{\SetFigFont{12}{14.4}{\familydefault}{\mddefault}{\updefault}{\color[rgb]{0,0,0}$j$}%
}}}
\put(5146,-4006){\makebox(0,0)[lb]{\smash{\SetFigFont{12}{14.4}{\familydefault}{\mddefault}{\updefault}{\color[rgb]{0,0,0}$k$}%
}}}
\put(5161,-6571){\makebox(0,0)[lb]{\smash{\SetFigFont{12}{14.4}{\familydefault}{\mddefault}{\updefault}{\color[rgb]{0,0,0}$k$}%
}}}
\put(4336,-5536){\makebox(0,0)[lb]{\smash{\SetFigFont{29}{34.8}{\familydefault}{\mddefault}{\updefault}{\color[rgb]{0,0,0}$=$}%
}}}
\end{picture}
 
\caption{ $\underline{(i,j)}.\underline{(j,k)}=\underline{(j,k)}.\underline{(i,k)}$
\label{relation standard}}
\end{center}
\end{figure}
\item[$\bullet$]
Si $r=(x,\xi)$ et $t=(x-n,\xi)$ avec $x\in X$ et $\xi\in \Xi$,
la relation est
$\underline{(x,\xi)}.\underline{(x-n,\xi)}=
\underline{(x+n,\xi)}.\underline{(x,\xi)}$.
La relation correspondante est vraie dans $\Bn$:
on consid\`ere l'automorphisme de $M(P_c)$ induit par l'identit\'e sur $X$ et la
translation de $n$ sur $\Xi$.
La relation revient \`a dire que
$\underline{(x,\xi)}.\underline{(x-n,\xi)}$
est invariant par cet automorphisme.
Dans $\Bn$ on consid\`ere l'automorphisme 
induit par une isotopie qui est
l'identit\'e sur le bord ext\'erieur de la couronne et fait
tourner l'autre bord de la couronne d'un tour dans le sens positif.
L'application de l'\'enonc\'e est compatible avec ces automorphismes.
On peut supposer $x\in [1,n]$;
en appliquant une puissance convenable de ces deux automorphismes on peut
ramener $\xi$ dans $[1,n]$;
la relation r\'esulte alors de la figure \ref{relation speciale}
qui permet de voir l'invariance
cherch\'ee.

\begin{figure}[htbp]
\begin{center}
 
\begin{picture}(0,0)%
\includegraphics{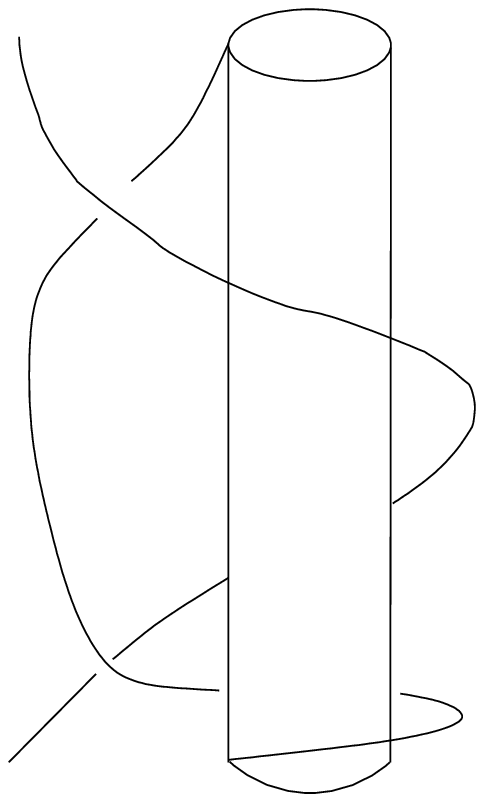}%
\end{picture}%
\setlength{\unitlength}{3947sp}%
\begingroup\makeatletter\ifx\SetFigFont\undefined%
\gdef\SetFigFont#1#2#3#4#5{%
  \reset@font\fontsize{#1}{#2pt}%
  \fontfamily{#3}\fontseries{#4}\fontshape{#5}%
  \selectfont}%
\fi\endgroup%
\begin{picture}(2330,3777)(1326,-7267)
\put(2251,-7111){\makebox(0,0)[lb]{\smash{\SetFigFont{12}{14.4}{\familydefault}{\mddefault}{\updefault}{\color[rgb]{0,0,0}$\xi$}%
}}}
\put(2251,-3661){\makebox(0,0)[lb]{\smash{\SetFigFont{12}{14.4}{\familydefault}{\mddefault}{\updefault}{\color[rgb]{0,0,0}$\xi$}%
}}}
\put(1331,-3646){\makebox(0,0)[lb]{\smash{\SetFigFont{12}{14.4}{\familydefault}{\mddefault}{\updefault}{\color[rgb]{0,0,0}$x$}%
}}}
\put(1326,-7131){\makebox(0,0)[lb]{\smash{\SetFigFont{12}{14.4}{\familydefault}{\mddefault}{\updefault}{\color[rgb]{0,0,0}$x$}%
}}}
\end{picture}
 
\caption{ $\underline{(x,\xi)}.\underline{(x,\xi+n)}$\label{relation speciale}}
\end{center}
\end{figure}

\end{enumerate}
\end{proof}
Les propositions \ref{B to G} et \ref{G to B} prouvent le th\'eor\`eme
\ref{G isom B}

Nous avons ainsi d\'efini pour chaque choix d'un \'el\'ement de Coxeter $c$ un
sous-mono{\"\i}de $M(P_c)$ de $\Bn$.
Deux tels mono{\"\i}des correspondant \`a des ensembles 
$\{a,\ldots,l\}$ et $\{\alpha,\ldots,\lambda\}$ de m\^emes cardinaux respectifs
sont isomorphes. L'\'echange de $X$ et de $\Xi$ est aussi un isomorphisme (qui
se traduit dans la repr\'esentation g\'eom\'etrique
par une rotation de $\pi$ de la
bande). On en d\'eduit que la classe d'isomorphisme
de $M(P_c)$ ne d\'epend que de
la partition de $n$ en $|X\cap [1,n]|+|\Xi\cap[1,n]|$. R\'eciproquement:
\begin{theoreme} L'application qui \`a $c$ associe
l'ensemble $\{|X\cap [1,n]|,|\Xi\cap[1,n]|\}$
induit une bijection des classes d'isomorphisme des mono{\"\i}des
$M(P_c)$ de la section \ref{monoide dual}
sur les partitions de $n$ en deux parties.
\end{theoreme}
Ce th\'eor\`eme r\'esulte des consid\'erations
qui pr\'ec\`edent et de la proposition
suivante qui montre que $|X\cap[1,n]|$ et $|\Xi\cap [1,n]|$ sont d\'etermin\'es
par l'action de la conjugaison par $\uc$ sur les atomes.
\begin{proposition}\label{orbites de c}
L'orbite par la conjugaison par $\uc$
d'un atome $\underline{(i,j)}$ est finie de cardinal $|X\cap[1,n]|$ (resp.
$|\Xi\cap[1,n]|$) si $i$ et $j$ sont tous
deux dans $X$ (resp. tous deux dans $\Xi$) et infinie si $i\in X$ et $j\in
\Xi$.
\end{proposition}
\begin{proof} Par \ref{conjugaison par c}, si $r$ et $r'$ sont des r\'eflexions
de $P_c$ on a $\ur\uc=\uc\ur'$ si et seulement si $rc=cr'$.
On a vu dans la d\'emonstration de \ref{presentation} que la conjugaison par $c$
revient \`a d\'ecaler $X$ d'une position dans le sens croissant et
$\Xi$ d'une position dans le sens
d\'ecroissant.
On en d\'eduit le r\'esultat.
\end{proof}

\section{Une structure \`a la Garside pour les groupes d'Artin-Tits affines de type $\tilde A$}\label{garside}

Gardons les notations des sections pr\'ec\'edentes.
Nous allons \'etudier les propri\'et\'es de la divisibilit\'e dans $P_c$ dans
le cas o\`u $c=s_1s_2\ldots s_n$. Dans ce cas on a
$c=(2,3,\ldots,n)_{[1]}(1)_{[-1]}$ et les ensembles $\Xi$ et $X$
sont respectivement $\Xi=\{z\in\bbZ\mid z\equiv1\pmod n\}$ et $X=\{z\in\bbZ\mid z\not\equiv1\pmod n\}$.

La propri\'et\'e fondamentale du mono{\"\i}de $M(P_c)$ dans ce cas
(th\'eor\`eme \ref{M(P_c) quasi-Garside})
est une cons\'equence de la proposition suivante:
\begin{proposition}\label{ppcm d'atomes} Si $c=s_1\ldots s_n$,
deux atomes quelconques de $P_c$ ont un ppcm dans $P_c$.
\end{proposition}
\begin{proof}
Consid\'erons deux atomes
$r=(x,y)$ et $r'=(x',y')$.
On sait que $r$ (resp. $r'$) divise un \'el\'ement
$p\in P_c$ si et seulement si $x$ et $y$ (resp. 
$x'$ et $y'$) sont dans la m\^eme partie de la partition associ\'ee \`a $p$
comme dans \ref{partition}.
Il y a deux cas.
\begin{itemize}
\item[$\bullet$]
Si $r$ et $r'$ sont sans croisement, ils
commutent et $rr'$ est dans $P_c$ par \ref{diviseurs de c}.
Dans la d\'ecomposition d'un multiple commun
$m$ de $r$ et $r'$ en produit de diviseurs
\'el\'ementaires, ou bien $r$
et $r'$ divisent deux diviseurs \'el\'ementaires  de $m$
diff\'erents et $rr'$ divise alors
$m$, ou bien $r$ et $r'$
divisent le m\^eme diviseur \'el\'ementaire $c'$ de $m$ qui est
un \'el\'ement de Coxeter d'un sous-groupe \quasiparab\  $W'$ et on conclut
que $rr'$ divise $c'$ soit par
\ref{diviseurs de c} si $W'$ est de type $\tilde A$ soit par
\cite[1.8]{BDM} si $W'$ est de type $A$. Donc le ppcm de $r$ et $r'$ existe et
vaut $rr'$.
\item[$\bullet$]
Si $r$ et $r'$ se croisent (\ie\ ne sont pas sans croisement), supposons
$x<y$, $x'<y'$ et $x$ et $x'$ dans $[1,n]$; si $m$ est un
multiple commun de $r$ et $r'$, 
alors $\{x,y\}$  et $\{x',y'\}$ ne peuvent pas \^etre dans deux parties
distinctes de la partition associ\'ee \`a $m$ car ces deux parties ne seraient
pas sans croisement.
Donc il existe un diviseur \'el\'ementaire $m'$ qui divise $m$ et qui est
multiple de $r$ et $r'$. Distinguons trois cas.

Si l'ensemble $\{x,y,x',y'\}\cap X$
est inclus dans un intervalle de
longueur strictement inf\'erieure \`a $n$, la partition dont les seules
parties non triviales sont $\{x,y,x',y'\}+kn$ avec $k\in\bbZ$
d\'efinit un \'el\'ement de Coxeter $c'$
d'un \quasiparab\ de type $A$.

Si  l'ensemble $\{x,y,x',y'\}\cap X$ n'est pas inclus dans un intervalle de longueur strictement inf\'erieure \`a $n$ et
si $\{x,y,x',y'\}\cap \Xi$ est non vide, la partition dont la seule partie non
triviale est $\{x,y,x',y'\}+n\bbZ$ d\'efinit un \'el\'ement
de Coxeter $c'$ d'un \quasiparab\ de type $\tilde A$.

Si  l'ensemble $\{x,y,x',y'\}\cap X$ n'est pas inclus dans un intervalle de longueur strictement inf\'erieure \`a $n$ et
si $\{x,y,x',y'\}\cap \Xi$ est vide, notons $c'$
l'\'el\'ement de Coxeter d'un \quasiparab\ de type $\tilde A$ d\'efini par la
partition dont la seule partie non triviale est
$\{1,x,y,x',y'\}+n\bbZ$ et remarquons que la partie associ\'ee \`a $m'$ \'etant
infinie doit contenir $1$.

Dans les trois cas l'\'el\'ement $c'$
divise $m'$ par \ref{diviseurs de c}  appliqu\'e \`a $m'$
et il est multiple de
$r$ et $r'$. C'est donc le ppcm de $r$ et $r'$ dans $P_c$.
\end{itemize}
\end{proof}
On peut alors appliquer la g\'en\'eralisation imm\'ediate suivante de \cite[0.5.2]{B}:
\begin{theoreme}\label{treillis} Soit $(G,R)$ un groupe positivement engendr\'e, soit
$c$ un \'el\'ement \'equilibr\'e et soit $M(P_c)$  comme pr\'ec\'edemment.
Supposons que deux \'el\'ements quelconques de $R$ ont un ppcm dans $P_c$; alors la
divisibilit\'e \`a gauche et la divisibilit\'e \`a droite donnent \`a $M(P_c)$ 
deux structures de treillis.
\end{theoreme}
Dans ce contexte on peut v\'erifier que les r\'esultats de \cite[section 2]{BDM} et
\cite[0.5]{B} s'appliquent.

Nous appellerons structure quasi-Garside une
structure de mono{\"\i}de v\'eri\-fiant tous les axiomes de \cite[0.5.1]{B} sauf
la finitude du nombre d'atomes.
Nos axiomes seront donc:
\begin{definition}\label{quasi-garside}
Un mono{\"\i}de $M$ est dit quasi-Garside si
\begin{enumerate}
\item Pour tout $m\in M$ le nombre de facteurs dans un
produit \'egal \`a $m$ est born\'e.
\item $M$ est simplifiable \`a gauche et \`a droite.
\item La divisibilit\'e \`a gauche et la divisibilit\'e \`a droite donnent \`a $M$ deux
structures de treillis.
\item Il existe un \'el\'ement $\Delta$ (\'el\'ement de Garside) dont l'ensemble des diviseurs \`a
gauche est \'egal \`a l'ensemble des diviseurs \`a droite et engendre $M$.
\end{enumerate}
\end{definition}
L'axiome (i) signifie que le mono{\"\i}de est atomique au sens par exemple de \cite[0.2.2]{B};
il est \'equivalent \`a l'existence
d'une longueur $l$ sur le mono{\"\i}de telle que $l(ab)\geq l(a)+l(b)$ pour tout couple
$(a,b)$.

Les axiomes (i), (ii) et (iv) \'etant v\'erifi\'es par un mono{\"\i}de d\'efini comme plus haut 
\`a partir d'un groupe engendr\'e et d'un \'el\'ement \'equilibr\'e,
l'\'enonc\'e \ref{treillis} devient:
\begin{theoreme}\label{M(P_c) quasi-Garside}
Sous les hypoth\`eses de \ref{treillis},
le mono{\"\i}de $M(P_c)$ est un mono\-{\"\i}de quasi-Garside avec $\uc$ comme
\'el\'ement de Garside.
\end{theoreme}
Ce th\'eor\`eme s'applique en particulier \`a un \'el\'ement $c$ comme dans
\ref{ppcm d'atomes}. Un tel choix d\'efinit donc une structure quasi-Garside
sur le groupe d'Artin-Tits de type $\tilde A_{n-1}$.

Le r\'esultat suivant montre que le choix de $c$ fait dans \ref{ppcm d'atomes}
est \`a isomorphisme pr\`es
le seul pour lequel la divisibilit\'e a une structure de treillis.
\begin{proposition}
Si $c$ est un \'el\'ement de Coxeter d'un groupe de type $\tilde A_{n-1}$, le
mono{\"\i}de $M(P_c)$ muni de l'ordre de la divisibilit\'e a une structure de
treillis (et est donc un mono{\"\i}de quasi-Garside)
si et seulement si l'un des deux ensembles $\Xi$ ou $X$ est r\'eduit \`a
un seul \'el\'ement modulo $n$.
\end{proposition}
\begin{proof}L'\'echange de $\Xi$ et $X$ d\'efinit un isomorphisme des mono{\"\i}des
correspondants. D'autre part
deux \'el\'ements de Coxeter tels que les ensembles $X$ correspondants
ont m\^eme nombre d'\'el\'ements modulo $n$ sont conjugu\'es, donc dans ce cas aussi
les mono{\"\i}des sont isomorphes. On en d\'eduit par \ref{ppcm d'atomes} que si
$\Xi$ ou $X$ a un seul \'el\'ement modulo $n$, on a bien une structure de
treillis. Inversement, supposons que modulo $n$, \`a la fois $\Xi$ et $X$ ont au
moins deux \'el\'ements. Soient $a<b$ dans $X\cap[1,n]$ (resp. $\alpha\neq \beta$ dans
$\Xi\cap[1,n]$). Les \'el\'ements $(a,b)$ et $(b,a+n)$
divisent $(a,b)_{[1]}(\alpha)_{[-1]}$ et $(a,b)_{[1]}(\beta)_{[-1]}$ qui sont
de longueur 3 et n'ont aucun diviseur commun de longueur 2 d'apr\`es \ref{diviseurs de c}.
Donc $(a,b)$ et $(b,a+n)$ n'ont pas de ppcm.
\end{proof}
Remarquons que dans l'exemple pr\'ec\'edent c'est le
dernier cas de la preuve de \ref{ppcm d'atomes}, dans lequel on a d\^u
introduire le \quasiparab\ de type $\tilde A$ d\'efini par
$\{1,x,y,x',y'\}+n\bbZ$,
qui est en d\'efaut. En fait les raisonnements des autres cas
s'appliquent pour tout \'el\'ement de Coxeter $c$ mais
ce dernier cas utilise le
fait que $\Xi$ est un singleton modulo $n$.

Donnons quelques
cons\'equences de l'existence d'une structure quasi-Gar\-side sur
$\Bn$. Ces cons\'equences sont de simples applications des propri\'et\'es g\'en\'erales
des mono{\"\i}des de Garside dont on v\'erifie qu'elles sont encore valables dans le cadre 
quasi-Garside.
Le premier est l'existence de formes normales telles que dues \`a Garside (voir par exemple
\cite{michel} ou \cite[section 2]{BDM}, voir aussi \cite{Charney}).
\begin{proposition}
\begin{enumerate}
\item
Tout \'el\'ement de $\Bn$ s'\'ecrit de fa\c con unique $a\inv b$ o\`u 
$a$ et $b$ sont des \'el\'ements de $M(P_c)$ premiers entre eux.
\item
Tout \'el\'ement de $M(P_c)$ s'\'ecrit de fa\c con unique $a_1a_2\ldots a_k$ o\`u pour
$i=1,\ldots,k$ l'\'el\'ement $a_i\in\uP_c$ est un (le) diviseur maximal dans
$\uP_c$ du produit $a_ia_{i+1}\ldots a_k$ et $a_k\neq 1$.
\end{enumerate}
\end{proposition}
Nous prouvons maintenant:
\begin{proposition}
Le centre de $\Bn$ est trivial.
\end{proposition}
\begin{proof}
La d\'emonstration suit les m\^emes id\'ees que la d\'e\-mons\-tration classique
pour les mono{\"\i}des de tresses ou que la d\'emonstration de \cite[4.1]{picantin}; ces
d\'emonstrations ne s'appliquent pas telles quelles car elles supposent qu'il y
a un nombre fini d'atomes. On utilise le lemme suivant:
\begin{lemme} Soit $M$ un mono{\"\i}de quasi-Garside et soit 
$b$ un \'el\'ement quasi-central de $M$, c'est-\`a-dire tel qu'il existe
un automorphisme $\tau$ de $M$ v\'erifiant $xb=b\tau(x)$ pour tout $x\in M$;
soit $x$ un diviseur \`a gauche de $b$ et $y\in M$; posons $\ppcm(x,y)=yz$ avec
$z\in M$: alors $z$ divise $b$ \`a gauche.
\end{lemme}
\begin{proof}
Comme $x$ divise $b$ il divise $b\tau(y)=yb$, donc $yz$ divise $yb$ et par
simplifiabilit\'e $z$ divise $b$.
\end{proof}
On en d\'eduit la proposition:
Soit $g\in\Bn$ central. On peut \'ecrire $g=\uc^nb$ avec $b\in M(P_c)$  non
divisible par $\uc$ et $n\in\bbZ$
convenable. On a alors $b$ quasi-central. Montrons par l'absurde que $b=1$.
Sinon, soit $r$ une r\'eflexion de $P_c$ telle que $\ur$ divise
$b$. Pour tout triplet de r\'eflexions $(r,s,t)$ correspondant \`a un triangle direct
comme dans la preuve de \ref{G to B} on a
$\ur.\us=\us.\ut=\ut.\ur=\ppcm(\ur,\us)=\ppcm(\us,\ut)$. On en d\'eduit par le lemme que
$\ur$ et $\ut$ divisent $b$. Comme \`a partir de $r$, de proche en proche on
peut faire appara{\^\i}tre n'importe quelle r\'eflexion $r'$ de $P_c$ dans un triangle
direct, on en d\'eduit que $b$ est multiple de tous les atomes de $M(P_c)$ donc
est multiple de $\uc$, ce qui est contradictoire. Tout \'el\'ement  central est donc
une puissance de $\uc$. Or la conjugaison par $\uc$ est un automorphisme d'ordre infini
(\cf\ \ref{presentation} ou \ref{orbites de c}),
donc aucune puissance de $\uc$ autre que $\uc^0$ n'est centrale.
\end{proof}

Donnons une derni\`ere cons\'equence de l'existence de la structure quasi-Garside.
\begin{proposition} \label{centralisateur}
Soient $c=s_1\ldots s_n$ et $\uc=\us_1\ldots \us_n$ comme pr\'ec\'edemment.
Le centralisateur de $\uc^h$ dans le groupe d'Artin-Tits de type
$\tilde A_{n-1}$ engendr\'e par $\us_1,\ldots,\us_n$ est un
groupe d'Artin-Tits de type $B_{\pgcd(h,n-1)}$.
\end{proposition}
\begin{proof}
La th\'eorie g\'en\'erale des mono{\"\i}des de (quasi-)Garside \'enonce que les points
fixes d'un automorphisme $\sigma$
dans un mono{\"\i}de $M$ de (quasi-)Garside forment un mono{\"\i}de $M^\sigma$ de
(quasi-)Garside, avec pour atomes certains des
ppcm des orbites des atomes et m\^eme \'el\'ement de Garside;
le groupe des fractions de ce mono\"\i de est le groupe des points fixes de $\sigma$ dans le
groupe des fractions de $M$ (\cf\ \cite[2.26]{BDM}
dont la d\'emonstration s'\'etend \`a la situation quasi-Garside).
Nous appliquons ceci \`a la conjugaison par $\uc^h$.
Le centralisateur de $\uc^h$ dans $G(P_c)$ est donc engendr\'e par les ppcm des
orbites des atomes sous la conjugaison par $\uc^h$.
Comme le montre le calcul fait dans \ref{presentation} la conjugaison par
$\uc$ envoie $\underline{(x,y)}$ sur $\underline{(x',y')}$ o\`u $x'$ et $y'$
s'obtiennent \`a partir de $x$ et de $y$ en translatant $X$ d'une position dans
le sens croissant et $\Xi$ d'une position dans le sens d\'ecroissant. Pour
simplifier au lieu d'indexer comme pr\'ec\'edemment les \'el\'ements de $X$ par les
entiers non congrus \`a 1 modulo $n$, nous renum\'erotons cons\'ecutivement les \'el\'ements de $X$,
en les notant $x_i$ avec $i\in \bbZ$, et nous faisons de m\^eme pour les
\'el\'ements de $\Xi$ qui seront not\'es $\xi_i$ avec $i\in\bbZ$. La translation de $n$
devient alors $x_i\mapsto x_{i+n-1}$ et $\xi_i\mapsto \xi_{i+1}$.
La conjugaison par $\uc^h$ envoie $\underline{(x_i,x_j)}$ sur
$\underline{(x_{i+h},x_{j+h})}$. L'orbite de $\underline{(x_i,x_j)}$
ne d\'epend donc que du pgcd de $h$ et de $n-1$.
La conjugaison par $\uc^h$ envoie $\underline{(\xi_0,x_j)}$  sur
$\underline{(\xi_0,x_{j+nh})}$. Le ppcm d'une telle orbite est
$\underline{(\xi_0)_{[-1]}(x_j,x_{j+k},x_{j+2k},\ldots,x_{j+n-1-k})_{[1]}}$
o\`u $k$ est le pgcd de $h$ et $n-1$. On voit que les ppcm des orbites d'atomes
ne d\'ependent que du pgcd de $n-1$ et $h$. Donc $\uc^h$ et $\uc^{\pgcd(n-1,h)}$
ont m\^eme centralisateur. On est donc ramen\'e au cas o\`u $h$ divise $n-1$. Le lemme suivant est le
cas particulier de la proposition quand $h=n-1$.
\begin{lemme}
Le centralisateur de $\uc^{n-1}$ est un groupe d'Artin-Tits de
type $B_{n-1}$ et l'image de $\uc$ dans ce groupe est un
\'el\'ement de Coxeter de ce groupe.
\end{lemme}
On a appel\'e \'el\'ement de Coxeter d'un groupe d'Artin-Tits le relev\'e canonique
d'un \'el\'ement de Coxeter du groupe de Coxeter.
\begin{proof}
Les \'el\'ements $\underline{(x_i,x_j)}$ sont centralis\'es par $\uc^{n-1}$. 
L'orbite de $\underline{(\xi_0,x_j)}$ se compose des \'el\'ements
$\underline{(\xi_0,x_{j+kn(n-1)})}$ avec $k\in\bbZ$. Le ppcm d'une telle
orbite est $\underline{(\xi_0)_{[-1]}(x_j)_{[1]}}$. Le centralisateur
$C(\uc^{n-1})$ de
$\uc^{n-1}$ est donc engendr\'e par ces \'el\'ements et toutes les relations
s'obtiennent en \'egalant les d\'ecompositions de $\uc$ comme produits de ces
g\'en\'erateurs. On a 
$$\uc=\underline{(x_1,x_2)}.\underline{(x_2,x_3)}\ldots\underline{(x_{n-2},x_{n-1})}.
\underline{(\xi_0)_{[-1]}(x_{n-1})_{[1]}}.$$

Pour construire un isomorphisme entre $C(\uc^{n-1})$ et le
groupe d'Artin-Tits $B(B_{n-1})$ de type $B_{n-1}$ nous revenons \`a l'interpr\'etation de
$B(\tilde A_{n-1})$ comme tresses dans une couronne. Remarquons que les
g\'en\'erateurs ci-dessus de $C(\uc^{n-1})$ sont des tresses telles que
le brin d'origine $\xi_0$ a pour extr\'emit\'e $\xi_0$ (voir figure \ref{figure3}). Or on peut interpr\'eter le
groupe de tresses de type $B_{n-1}$ comme le groupe de tresses \`a $n-1$
brins dans la couronne, les points de base \'etant les points
$x_1,x_2,\ldots,x_{n-1}$ (\cf, \cite{Lambropoulou}). On a donc un morphisme de $C(\uc^{n-1})$
dans $B(B_{n-1})$ par oubli du brin issu de $\xi_0$.
Les g\'en\'erateurs de $C(\uc^{n-1})$ s'envoient sur les g\'en\'erateurs de la
pr\'esentation duale du groupe $B(B_{n-1})$ et l'image de $\uc$ est l'\'el\'ement
de Coxeter de $B(B_{n-1})$. Donc les relations entre les g\'en\'erateurs de
ces deux groupes se correspondent. Le morphisme est donc un isomorphisme.
\end{proof}

\begin{figure}[htbp]
\begin{center}
 
\begin{picture}(0,0)%
\includegraphics{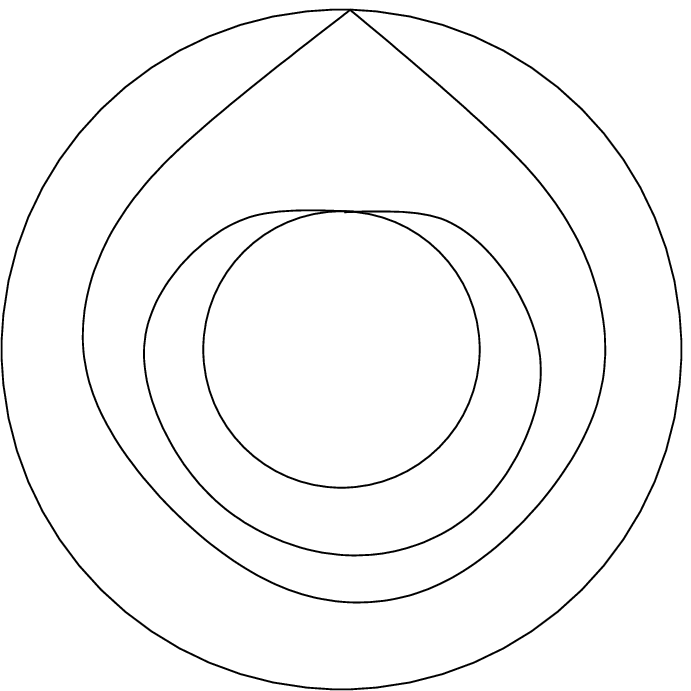}%
\end{picture}%
\setlength{\unitlength}{4144sp}%
\begingroup\makeatletter\ifx\SetFigFont\undefined%
\gdef\SetFigFont#1#2#3#4#5{%
  \reset@font\fontsize{#1}{#2pt}%
  \fontfamily{#3}\fontseries{#4}\fontshape{#5}%
  \selectfont}%
\fi\endgroup%
\begin{picture}(3124,3338)(3839,-6122)
\put(5317,-4107){\makebox(0,0)[lb]{\smash{\SetFigFont{12}{14.4}{\familydefault}{\mddefault}{\updefault}{\color[rgb]{0,0,0}$\xi_0$}%
}}}
\put(5213,-2941){\makebox(0,0)[lb]{\smash{\SetFigFont{12}{14.4}{\familydefault}{\mddefault}{\updefault}{\color[rgb]{0,0,0}$x_{n-1}$}%
}}}
\put(6890,-3964){\makebox(0,0)[lb]{\smash{\SetFigFont{12}{14.4}{\familydefault}{\mddefault}{\updefault}{\color[rgb]{0,0,0}$x_2$}%
}}}
\put(6405,-3299){\makebox(0,0)[lb]{\smash{\SetFigFont{12}{14.4}{\familydefault}{\mddefault}{\updefault}{\color[rgb]{0,0,0}$x_1$}%
}}}
\put(6959,-4850){\makebox(0,0)[lb]{\smash{\SetFigFont{12}{14.4}{\familydefault}{\mddefault}{\updefault}{\color[rgb]{0,0,0}$x_3$}%
}}}
\put(6681,-5595){\makebox(0,0)[lb]{\smash{\SetFigFont{12}{14.4}{\familydefault}{\mddefault}{\updefault}{\color[rgb]{0,0,0}$x_4$}%
}}}
\put(4029,-3333){\makebox(0,0)[lb]{\smash{\SetFigFont{12}{14.4}{\familydefault}{\mddefault}{\updefault}{\color[rgb]{0,0,0}$x_{n-2}$}%
}}}
\end{picture}
 
\caption{ $\underline{(\xi_0)_{[-1]}(x_{n-1})_{[1]}}$ \label{figure3}}
\end{center}
\end{figure}

Prouvons alors la proposition. Si $h$ divise $n-1$,
le centralisateur de $\uc^h$ est \'egal au centralisateur de $\uc^h$ dans le centralisateur de
$\uc^{n-1}$.
Le lemme permet donc de terminer la d\'emonstration de la proposition car le
centralisateur d'une puissance d'un \'el\'ement de Coxeter dans un groupe de tresses de type $B$
est connu par les r\'esultats de \cite{BDM}. 
\end{proof}
\begin{remarque}
On obtient des g\'en\'erateurs standards (\ie v\'erifiant les relations de tresses de type $B$)
du centralisateur de $\uc^{n-1}$ en prenant
$\underline{(x_1,x_2)}$, $\underline{(x_2,x_3)}$, \dots, $\underline{(x_{n-2},x_{n-1})}$, et
$\underline{(\xi_0)_{[-1]}(x_{n-1})_{[1]}}$.
D'autre part il est facile de voir que le
centralisateur de $(s_1s_2\ldots s_n)^{n-1}$ dans le groupe de Coxeter $W(\tilde
A_{n-1})$ est le groupe le groupe engendr\'e par
$(x_1,x_2)$, $(x_2,x_3)$, \dots, $(x_{n-2},x_{n-1})$, et
$(\xi_0)_{[-1]}(x_{n-1})_{[1]}$ donc est l'image du centralisateur de
$\uc^{n-1}$ dans le groupe d'Artin-Tits. La pr\'esentation de cette image s'obtient en
ajoutant aux relations de tresses de type $B$ le fait que les $n-2$ premiers
g\'en\'erateurs sont d'ordre 2 (le dernier est d'ordre infini).
\end{remarque}

\section*{Remerciements}
Je remercie tout particuli\`erement Eddy Godelle
pour ses nombreuses remarques pertinentes sur une version ant\'erieure
de cet article. Je remercie \'egalement Jean Michel
avec qui j'ai eu plusieurs discussions sur ces r\'esultats.

\end{document}